\documentclass[leqno,11pt]{amsart}

\usepackage{amssymb}
\usepackage{amsmath}
\usepackage{enumerate}
\usepackage[pdftex]{graphicx}
\usepackage{graphicx, color}

\def\R{\mathbb{R}}

\def\Q{\mathbb{Q}}

\def\lo{\mathbb{L}}

\newcommand{\arcsinh}{\mathop{\rm arcsinh }\nolimits}
\DeclareMathOperator{\arctanh}{arctanh}
\DeclareMathOperator{\arccoth}{arccoth}
\DeclareMathOperator{\arccosh}{arccosh}

\newtheorem{remark}{Remark}[section]
\newtheorem{theorem}{Theorem}[section]

\newtheorem{corollary}{Corollary}[section]

\newtheorem{example}{Example}[section]

\numberwithin{equation}{section}

\begin{document}

\title[Curves in $\lo^2$]{Curves in Lorentz-Minkowski plane with curvature depending on their position}

\author[I. Castro]{Ildefonso Castro}
\address{Departamento de Matem\'{a}ticas \\
Universidad de Ja\'{e}n, 23071 Ja\'{e}n, Spain \\
and Instituto de Matemáticas (IEMath-GR)
} 
\email{icastro@ujaen.es}

\author[I. Castro-Infantes]{Ildefonso Castro-Infantes}
\address{Departamento de Geometr\'{\i}a y Topolog\'{\i}a \\
Universidad de Granada, 
18071 Granada, Spain \\
and Instituto de Matemáticas (IEMath-GR)} 
\email{icastroinfantes@ugr.es}

\author[J. Castro-Infantes]{Jes\'{u}s Castro-Infantes}
\address{Departamento de Geometr\'{\i}a y Topolog\'{\i}a \\
Universidad de Granada \\
18071 Granada, Spain} \email{jcastroinf@correo.ugr.es}

\thanks{Research of the two first named authors was partially supported by 
a MEC-FEDER grant MTM2017-89677-P. Research of the third named author was partially supported by 
a MECD grant FPU16/03096.}

\subjclass[2010]{Primary 53A35; Secondary 14H50, 53B30, 74H05}

\keywords{Lorentzian curves, curvature, Singer's problem, sinusoidal spirals, elastic curves, grim-reaper curves.}

\date{}

\begin{abstract}\textit{}
Motivated by the classical Euler elastic curves, David A. Singer posed in 1999 the problem of determining a plane curve whose curvature is given in terms of its position.
We propound the same question in Lorentz-Minkowski plane, focusing on spacelike and timelike curves.
In this article, we study those curves in $\lo^2$ whose curvature depends on the Lorentzian pseudodistance from the origin, and those ones whose curvature depends on the Lorentzian pseudodistance through the horizontal or vertical geodesic to a fixed lightlike geodesic. Making use of the notions of {\em geometric angular  momentum} (with respect to the origin) and {\em geometric linear  momentum} (with respect to the fixed lightlike geodesic) respectively, we get 
two abstract integrability results to determine such curves through quadratures. In this way, we find out several new families of Lorentzian spiral, special elastic and grim-reaper curves whose intrinsic equations are expressed in terms of elementary functions. In addition, we provide uniqueness results for the generatrix curve of the Enneper's surface of second kind and for Lorentzian versions of some well known curves in the Euclidean setting, like the Bernoulli lemniscate, the cardioid, the sinusoidal spirals and some non-degenerate conics.  We are able to get arc-length parametrizations of them and they are depicted graphically.
\end{abstract}

\maketitle

\section{Introduction}
David A. Singer proposed in 1999 (cf.\ \cite{S99}) the following problem: 
\begin{quote}
{\em Can a plane curve be determined if its curvature is given in terms of its position?}
\end{quote}
This question was motivated by the fact that for the classical Euler elastic curves (cf.~\cite{S08} for instance),
their curvature is proportional to one of the coordinate functions.
Singer proved (see Theorem 3.1 in \cite{S99}) that the
problem of determining a curve whose curvature is $\kappa (r)$, where $r$ is the distance from the origin, is solvable by quadratures when $r\kappa (r)$ is a continuous function.
Singer himself studied in \cite{S99} the simple case $\kappa (r)=r$, where elliptic integrals appear in, illustrating the fact that the corresponding differential equations are integrable by quadratures does not mean that it is easy to perform the integrations. In that paper, only at the very pleasant special case of the classical Bernoulli lemniscate, the aforementioned integrals became elementary. Nevertheless, the curvature of many famous plane curves, such as conic sections, Cassini ovals, Sturmian and sinusoidal spirals, depends only on distance $r$ from a given point (that we fix as the origin).

The recent literature devoted to the study of particular cases of Singer's posed problem consisting on determining plane curves $\alpha\!=\!(x,y)$ given $\kappa \!=\! \kappa(r)$, $r=\sqrt{x^2+y^2}$, includes several papers like \cite{DVM09}, \cite{VDM09}, \cite{MHDV10}, \cite{MHDV11a}, \cite{MHDV11b}, \cite{MHM14} or \cite{MMH14}. In addition, the authors have studied the cases $\kappa \!=\! \kappa(y)$ and $\kappa \!=\! \kappa(r)$ in \cite{CCI16} and \cite{CCIs17} respectively, for a large number of prescribed curvature functions. 

Our purpose in this paper is to propound the study of Singer's problem in Lorentz-Minkowski plane;
that is, to try to determine those curves $\gamma =(x,y)$ in $\lo^2$ whose curvature $\kappa$ depends on some given function $\kappa =\kappa(x,y)$. We focus on spacelike and timelike curves, since the curvature $\kappa $ is in general not well defined on lightlike points, and because lightlike curves in $\lo^2$ are segments parallel to the straight lines determining the light cone.
When the ambient space is $\lo^2$, it is still valid the fundamental existence and uniqueness theorem. It states that a spacelike or timelike curve is uniquely determined (up to Lorentzian transformations) by its curvature given as a function of its arc-length, although in practice it is very difficult to find the curve explicitly at most cases, such as it happens in the Euclidean case. But now, referred to the Lorentzian version of Singer's problem, our knowledge is much more restricted in comparison with the Euclidean case. In fact, we can only mention the articles \cite{IUM15} and \cite{UIM16} in this line, both devoted to Sturmian spiral curves. The authors initiated in \cite{CCIs18} the study of the spacelike and timelike curves in $\lo^2$ satisfying $\kappa =\kappa(y)$ or $\kappa =\kappa(x)$, both conditions geometrically interpreted as that the curvature is expressed in terms of the Lorentzian pseudodistance to timelike or spacelike fixed geodesics. 
In this way, we get in \cite{CCIs18} a complete description of almost all the elastic curves in $\lo^2$ and 
provide the Lorentzian versions of catenaries and some grim-reaper curves.

From a geometric-analytic point of view, we afford in this article the following two cases of Singer's problem in the Lorentzian setting: 
For a unit-speed parametrization of a regular curve $\gamma =(x,y)$ in $\lo^2:=(\R^2,g=-dx^2+dy^2)$, 
we prescribe the curvature with the geometric extrinsic conditions $\kappa\!=\!\kappa (\rho)$ and $\kappa\!=\!\kappa (v)$, where $\rho: =\! \sqrt{|g(\gamma, \gamma)|} \!=\! \sqrt{|-x^2+y^2|}\geq 0$ and $v\!=\!y-x$. 
The first condition is interpreted geometrically as 
those curves whose curvature depends on the Lorentzian pseudodistance from the origin, and the second one as those ones whose curvature depends on the Lorentzian pseudodistance through the horizontal or vertical geodesic to a fixed lightlike geodesic (see Section 2 for details).
We aim to determine explicitly the analytic representation of the arc-length parametrization $\gamma (s) $ and, as a consequence, its intrinsic equation $\kappa = \kappa (s)$.

Singer's proof of the aforementioned Theorem 3.1 in \cite{S99} is based in giving to such  a curvature $\kappa = \kappa (r)$ an interpretation of a central potential in the plane and finding the trajectories by the standard methods in classical mechanics. On the other hand, since the curvature $\kappa$ may be also interpreted as the tension that $\gamma $ receives at each point as a consequence of the way it is  immersed in $\lo^2$, we make use of the notions of {\em geometric angular  momentum} (with respect to the origin) of $\gamma$ when $\kappa\!=\!\kappa (\rho)$ and {\em geometric linear  momentum} (with respect to the fixed lightlike geodesic) of $\gamma$ when $\kappa\!=\!\kappa (v)$  in order to get 
two  abstract integrability results (Theorems~\ref{quadratures rho} and \ref{quadratures v}) in the same spirit of Theorem 3.1 in \cite{S99}. We show that the problem of determining such a curve is solvable by quadratures if $\rho \,\kappa (\rho)$ and $\kappa (v)$ are continuous functions at each case respectively. In addition, both geometric momenta determine uniquely such a curve (up to certain Lorentz transformations). In general, one finds in most cases great difficulties (see for example Remark~\ref{difficulties r}) to carry out the computations. Hence we focus on finding interesting curves for which the required computations can be achieved {\em explicitly}, in terms of elementary functions, in a straightforward way from our two main results. We pay attention to identify, compute and plot such examples. 

In this way, we are successful in Section 4 by studying the most regular curves in the Lorentz-Minkowski plane satisfying the geometric condition $\kappa (\rho)=2\lambda + \mu/\rho$, $\lambda, \mu \in \R$, generalizing the Sturmian spirals studied in \cite{IUM15} corresponding to $\lambda =0$. Moreover, we find out the Lorentzian versions of some classical interesting curves in the Euclidean context corresponding to the sinusoidal spirals, characterized in Corollary \ref{cor:sinusoidal} by their geometric angular momentum $\mathcal K(\rho)=  \frac{\lambda}{n+1}\rho^{n+1} $ (and curvature  $\kappa (\rho)=\lambda \, \rho^{n-1}),\, \lambda >0$, $ n\neq 0, \, n\neq -1$.
They include for particular values of $n$ the Lorentzian counterparts of the Bernoulli lemniscate, the cardioid and some non-degenerate conics (Remark \ref{re:conics}). 

On the other hand, we completely describe in Section 6 the spacelike and timelike curves in $\lo^2$ satisfying
$\kappa(v)  = a v + b, \ a \neq 0, \, b \in \R$, recovering in this new way the special elastic curves described in Section 3 of \cite{CCIs18} with energy $E=\sigma ^2 /4$, being $\sigma $ the tension of the elastica. 
Moreover, we classify in Section 7 the causal curves in $\lo^2$ satisfying
$\kappa(v)  = a / v^2,  a > 0$, and provide a uniqueness result (Corollary \ref{cor:Ennep2}) for the generatrix curve of the Enneper's surface of second kind described in \cite{Ko93}. Finally, in Section 8 we determine those spacelike and timelike curves in $\lo^2$ satisfying $\kappa(v)  = a\,  e^v, a > 0$, characterizing in Corollary \ref{cor:grims} new Lorentzian grim-reapers, i.e.\ Lorentzian curves that satisfy certain translating-type soliton equations.

\vspace{0.3cm}

\section{Spacelike and timelike curves in Lorentz-Minkowski plane}
We denote $\lo^2\!:=\!(\R^2, g\!=\!-dx^2+dy^2)$ the Lorentz-Minkowski plane, where $(x,y)$ are the rectangular coordinates on $\lo^2$. We say that a non-zero vector $v\in \l ^2$ is spacelike if $g(v,v)>0$, lightlike if $g(v,v)=0$, and timelike if $g(v,v)<0$.

Let $\gamma \!=\!(x,y) \!:\! I\subseteq\R \rightarrow \R^2$ be a curve. We say that $\gamma=\gamma(t)$ is spacelike (resp.\ timelike) if the tangent vector $\gamma'(t)$ is spacelike (resp.\ timelike) for all $t\in I$. A point $\gamma (t)$ is called a lightlike point if $\gamma'(t)$ is a lightlike vector.
We study in this section geometric properties of curves that have no lightlike points, because the curvature is not in general well defined at the lightlike points. 

Hence, let $\gamma \!=\!(x,y)$ be a spacelike (resp.\ timelike) curve parametrized by arc-length; that is, $g(\dot \gamma(s),\dot \gamma(s))= 1$ (resp.\ $g(\dot \gamma(s),\dot \gamma(s))= -1$) $\forall s\in I $, where $I$ is some interval in $\R$. Here $\ \dot{}\ $ means derivation with respect to $s$. 
We will say in both cases that $\gamma=\gamma(s)$ is a unit-speed curve.

\begin{figure}[h!]
\begin{center}
\includegraphics[height=4cm]{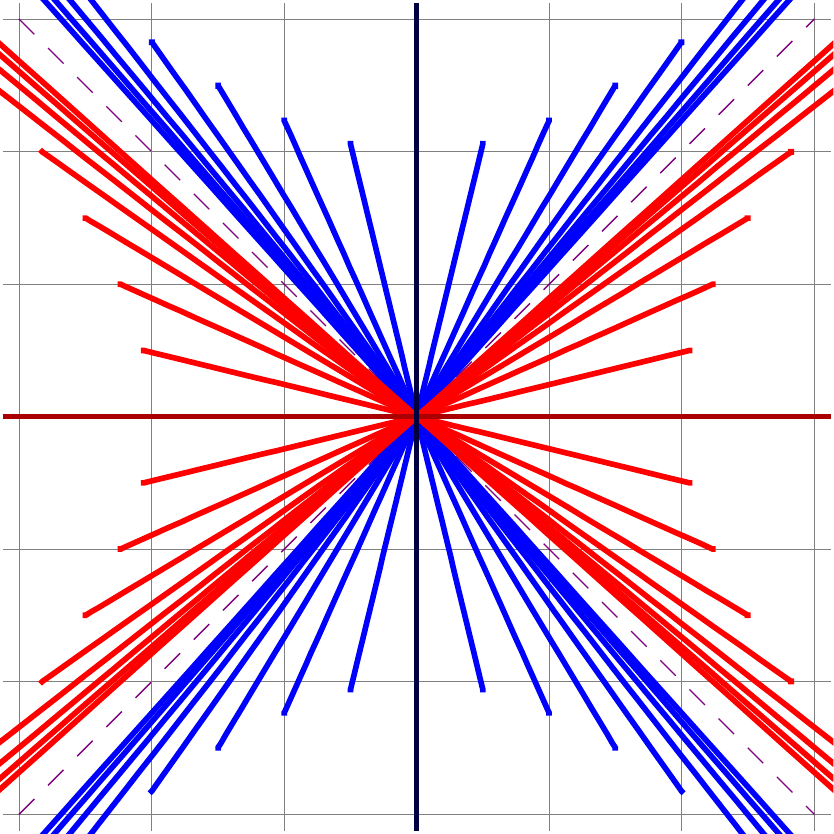}
\caption{Spacelike (blue) and timelike (red) geodesics in $\lo^2$.}
\label{Fig:geodesics}
\end{center}
\end{figure}

We define the Frenet dihedron in such a way that the curvature $\kappa$ has a sign and then it is only preserved by direct rigid motions (see \cite{Lo14}):
Let $T=\dot{\gamma}=(\dot x,\dot y)$ be the tangent vector to the curve $\gamma$ and let $N=(\dot y, \dot x)$ be the corresponding normal vector. We remark that $T$ and $N$ have different causal character. Let $g(T,T)=\epsilon$, with $\epsilon= 1$ if $\gamma$ is spacelike, and $\epsilon= -1$ if $\gamma$ is  timelike. Then $g(N,N)=-\epsilon$. The (signed) curvature of $\gamma$ is the function $\kappa=\kappa(s)$ such that
\begin{equation}\label{FrenetEq1}
\dot T(s)= \kappa (s)  N(s),
\end{equation}
where
\begin{equation}\label{curvature}
\kappa(s)=-\epsilon g(\dot T(s),N(s))=\epsilon (\ddot x \dot y - \dot x \ddot y).
\end{equation}
The Frenet equations of $\gamma$ are given by \eqref{FrenetEq1} and
\begin{equation}\label{FrenetEq2}
 \dot N(s)= \kappa (s) T(s).
\end{equation}
It is possible, as it happens in the Euclidean case, to obtain a parametrization by arc-length of the curve $\gamma$ in terms of integrals of the curvature. Concretely, any spacelike curve $\alpha (s)$ in $\lo^2$ can be represented (up to isometries) by
\begin{equation}\label{spacelike curve}
\alpha(s)=\left( \int \sinh \varphi (s) ds , \int \cosh \varphi(s) ds \right), \ \frac{d\varphi(s)}{ds}=\kappa (s),
\end{equation}
and any timelike curve $\beta (s)$ in $\lo^2$ can be represented (up to isometries) by
\begin{equation}\label{timelike curve}
\beta(s)=\left( \int \cosh \phi (s) ds , \int \sinh \phi(s) ds \right), \ \frac{d\phi(s)}{ds}=\kappa (s).
\end{equation}
For example, up to a translation, any spacelike geodesic can be written as
\begin{equation}\label{spacelike geodesic}
\alpha_{\varphi_0}(s)=( \sinh \varphi_0 \, s ,  \cosh \varphi_0 \, s), \, s\in \R, \ \varphi_0\in \R,
\end{equation}
and any timelike geodesic can be written as
\begin{equation}\label{timelike geodesic}
\beta_{\phi_0}(s)=( \cosh \phi_0 \, s ,  \sinh \phi_0 \, s), \, s\in \R, \ \phi_0\in \R.
\end{equation}
On the other hand, the transformation $R_\nu: \lo^2 \rightarrow \lo^2, \, \nu \in \R $, given by
\begin{equation}\label{Rnu}
R_\nu (x,y)=(\cosh \nu \,x + \sinh \nu \,y, \sinh \nu \,x + \cosh \nu \,y ) 
\end{equation}
is an isometry of $\lo^2$ that preserves the curvature of a curve $\gamma $ and satisfies
$$ R_\nu \circ \alpha_{\varphi_0} = \alpha_{\varphi_0+\nu}, \  R_\nu \circ \beta_{\phi_0}= \beta_{\phi_0+\nu}. $$
We will refer $R_\nu$, $\nu \in \R$, as a $\nu$-orthochrone Lorentz transformation (see \cite{Lo14}).
In this way, any spacelike geodesic is congruent to $\alpha_0$, i.e.\ the $y$-axis, and any timelike geodesic is congruent to $\beta_0$, i.e.\ the $x$-axis (see Figure \ref{Fig:geodesics}).

\subsection{Curves in $\lo^2$ such that $\kappa =\kappa(\rho)$}
We are first interested in spacelike or timelike unit-speed curves $\gamma =(x,y)$ satisfying the geometric condition $\kappa =\kappa (\rho)$, where $\rho$ is the {\em Lorentzian pseudodistance from the origin} defined by 
$$\rho:=\sqrt{|g(\gamma, \gamma)|}=\sqrt{|-x^2+y^2|}\geq 0.$$
For this aim, we use for $\gamma$ what  we can call {\em pseudopolar coordinates} 
$(\rho, \nu)$, $\rho \geq 0$, $ \nu \in \R$ being the {\em orthochrone angle}. 

Since $g(\gamma, \gamma)=-x^2+y^2 = \pm \rho^2$, we will write:
\begin{equation}\label{pseudopolar0}
\begin{array}{l}
\gamma^+ \equiv \left\{ 
\begin{array}{ll}
x=\rho \, \sinh \nu, \, y= \rho \, \cosh \nu, & {\rm if \ } -x^2+y^2 \geq 0, y \geq 0  \\ 
x=-\rho \, \sinh \nu, \, y= -\rho \, \cosh \nu, & {\rm if \ } -x^2+y^2 \geq 0, y \leq 0 
\end{array}
\right. \\ \\
\gamma^- \equiv \left\{ 
\begin{array}{ll}
x=\rho \, \cosh \nu, \, y= \rho \, \sinh \nu, & {\rm if \ } -x^2+y^2 \leq 0, y \geq 0  \\  
x=-\rho \, \cosh \nu, \, y= -\rho \, \sinh \nu, & {\rm if \ } -x^2+y^2 \leq 0, y \leq 0 
\end{array}
\right.
\end{array}
\end{equation}
In fact, it will be enough obviously to consider the first and third cases, since the map $(x,y) \rightarrow (-x,-y)$ is an isometry of $\lo^2$.

When $ x=\rho \sinh \nu, \, y= \rho \cosh \nu$, we have that
$g(\gamma, \gamma) = \rho^2 $, $ g(\dot \gamma, \dot \gamma)=\dot \rho^2-\rho^2 \dot \nu^2 $ and $g(\gamma, N)=\rho^2 \dot \nu$. 
At a point $\gamma (s)=(\rho (s) \sinh \nu (s),\rho(s) \cosh \nu (s))$ on the curve, the {\em geometric angular momentum} $\mathcal K$  is given by
\begin{equation}\label{momentum}
\mathcal K(s) = \rho (s)^2 \, \dot \nu (s).
\end{equation}
In physical terms, using Noether's Theorem, $\mathcal K$ may be interpreted as the angular momentum with respect to the origin of a particle of unit mass with unit-speed and trajectory $\gamma $. 
Using that $\gamma$ is unit-speed, that is, $\dot \rho^2-\rho^2 \dot \nu^2 = \epsilon $, and \eqref{momentum}, we easily obtain that
\begin{equation}\label{differentials}
ds = \frac{\rho \, d\rho}{\sqrt{\mathcal K ^2+\epsilon \rho^2}}, \quad d\nu =\frac{\mathcal K \,  ds}{\rho^2}.
\end{equation}
Thus, given $\mathcal K=\mathcal K (\rho)$ as an explicit function, looking at \eqref{differentials} one may attempt to compute $\rho(s)$ and $\nu (s) $ in three steps: integrate to get $s=s(\rho)$, invert to get $\rho=\rho(s)$ and integrate to get $\nu =\nu(s)$. 
In addition, in this case $\mathcal K= g(\gamma, N)$ taking into account \eqref{momentum}. Using \eqref{FrenetEq2} and $g(\gamma, \gamma) = \rho^2 $, we deduce that $\frac{d\mathcal K}{ds}=\kappa \dot \rho \rho$ and, since we are assuming that $\kappa =\kappa (\rho)$, we finally arrive at 
\begin{equation}\label{K primitive}
d \mathcal K  =  \rho \kappa(\rho)d\rho,
\end{equation}
that is, $\mathcal K= \mathcal K (\rho)$ can be interpreted as an anti-derivative of $\rho \kappa (\rho)$.

\begin{figure}[h!]
\begin{center}
\includegraphics[height=4cm]{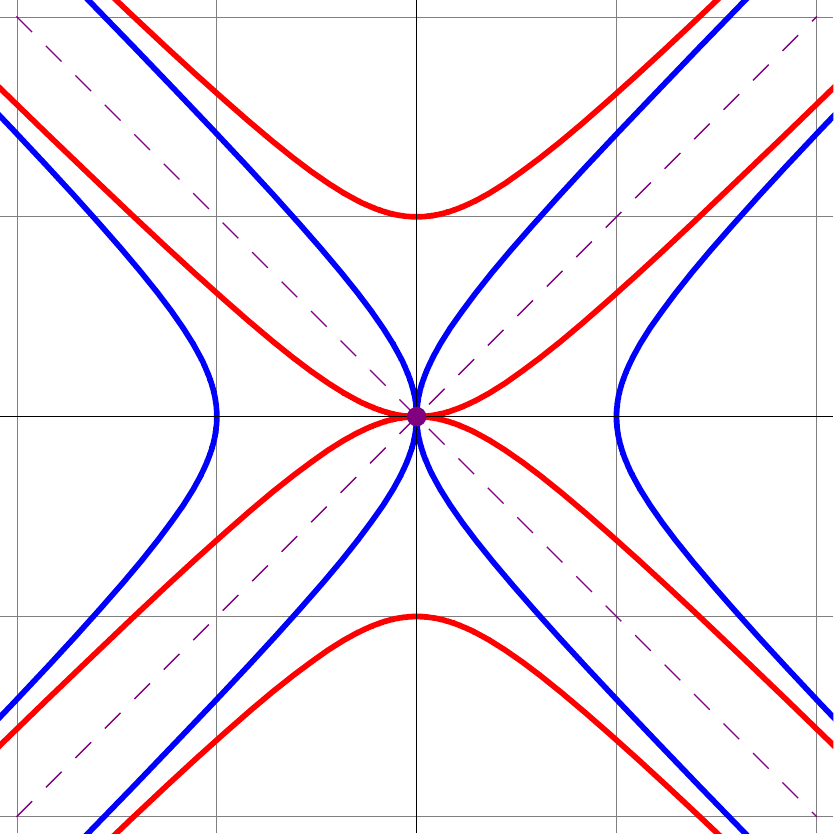}
\caption{Spacelike (blue) and timelike (red) pseudocircle with $\mathcal K (\rho)=\rho^2/2$ ($\kappa \equiv 1$)  in $\lo^2$.}
\label{Fig:circlesO}
\end{center}
\end{figure}

And when $ x=\rho \cosh \nu, \, y= \rho \sinh \nu$, since we then have that
$g(\gamma, \gamma) = -\rho^2 $, $ g(\dot \gamma, \dot \gamma)=-\dot \rho^2+\rho^2 \dot \nu^2 = \epsilon $ and $g(\gamma, N)=-\rho^2 \dot \nu$, with the same definition \eqref{momentum} for $\mathcal K$, \eqref{K primitive} holds and \eqref{differentials} changes into
\begin{equation}\label{differentials bis}
ds = \frac{\rho \, d\rho}{\sqrt{\mathcal K ^2-\epsilon \rho^2}}, \quad d\nu =\frac{\mathcal K \,  ds}{\rho^2}.
\end{equation}
As a summary, we have proved the following result in the spirit of Theorem 3.1 in \cite{S99}.
\begin{theorem}\label{quadratures rho}
Let $\kappa=\kappa(\rho)$ be a function such that $\rho \, \kappa (\rho)$ is continuous.
Then the problem of determining locally a spacelike or timelike curve in $\lo^2$
whose curvature is $\kappa(\rho )$ with geometric angular momentum $\mathcal K (\rho)$ satisfying \eqref{K primitive}
---$\rho$ being the (non constant) Lorentzian pseudodistance from the origin---
is solvable by quadratures considering the unit speed curve $\gamma$ in pseudopolar coordinates \eqref{pseudopolar0} given by
\begin{equation}\label{pseudopolar}
\begin{array}{c}
\gamma^+_\epsilon (s)=(\pm \rho^+_\epsilon (s) \sinh \nu^+_\epsilon (s),\pm \rho^+_\epsilon (s) \cosh \nu^+_\epsilon (s)), \\
\gamma^-_\epsilon(s)=(\pm \rho^-_\epsilon (s) \cosh \nu^-_\epsilon (s),\pm \rho^-_\epsilon (s) \sinh \nu^-_\epsilon (s)),
\end{array}
\end{equation}
where $\rho^+_\epsilon(s)$ and $\nu^+_\epsilon (s)$ 
(resp.\ $\rho^-_\epsilon(s)$ and $\nu^-_\epsilon (s)$)
are obtained through (\ref{differentials}) (resp.\ through (\ref{differentials bis})) after inverting $s=s(\rho)$,
being $\epsilon=1$ at the spacelike case and $\epsilon=-1$ at the timelike case.
Such a curve is uniquely determined by $\mathcal K (\rho)$ up to a $\nu$-orthochrone Lorentz transformation (and a translation of the arc parameter $s$).
\end{theorem}
\begin{remark}\label{c rho}
{\rm If we prescribe $\kappa=\kappa(\rho)$, the method described in Theorem~\ref{quadratures rho} clearly implies the computation of four quadratures, following the sequence:
\begin{enumerate}
\item[\rm (i)]  Anti-derivative of $\rho \, \kappa (\rho)$: $$\int \! \rho \, \kappa(\rho)d\rho= \mathcal K (\rho). $$
\item[\rm (ii-a)]  Arc-length parameter of $$\gamma^+_\epsilon (s)\!=\!(\pm \rho^+_\epsilon (s) \sinh \nu^+_\epsilon (s),\pm \rho^+_\epsilon (s) \cosh \nu^+_\epsilon (s)) $$ in terms of $\rho$: $$s=s(\rho)=\int \!\frac{\rho \, d\rho }{\sqrt{\mathcal K (\rho)^2+\epsilon \rho^2}},$$
where $\mathcal K (\rho)^2+\epsilon \rho^2 >0$, and inverting $s=s(\rho)$ to get $\rho=\rho^+_\epsilon(s)>0$.
\item[\rm (ii-b)]  Arc-length parameter of $$\gamma^-_\epsilon(s)\!=\!(\pm \rho^-_\epsilon (s) \cosh \nu^-_\epsilon (s),\pm \rho^-_\epsilon (s) \sinh \nu^-_\epsilon (s)) $$ in terms of $\rho$: $$s=s(\rho)=\int\!\frac{\rho \, d\rho }{\sqrt{\mathcal K (\rho)^2-\epsilon \rho^2}},$$
where $\mathcal K (\rho)^2-\epsilon \rho^2 >0$, and inverting $s=s(\rho)$ to get $\rho=\rho^-_\epsilon(s)>0$.
\item[\rm (iii)] Orthochrone angles $ \nu^+_\epsilon $ and $\nu^-_\epsilon$ in terms of $s$: 
$$\nu^+_\epsilon (s)=\int \frac{\mathcal K (\rho^+_\epsilon(s))}{\rho^+_\epsilon (s)^2}ds,
\quad
\nu^-_\epsilon (s)=\int \frac{\mathcal K (\rho^-_\epsilon(s))}{\rho^-_\epsilon (s)^2}ds,
$$
where $\rho^+_\epsilon (s)>0$ and $\rho^-_\epsilon (s)>0$ respectively.
\end{enumerate}
We note that we get a one-parameter family of curves in $ \lo^2$ satisfying $\kappa=\kappa(\rho)$ according to the geometric angular momentum chosen in (i). It will distinguish geometrically the curves inside a same family by their relative position with respect to the origin. We remark that we can recover $\kappa$ from $\mathcal K$ by means of $\kappa (\rho)= (1/\rho) d\mathcal K/d\rho$. In addition, parts (ii-a) and (ii-b) are consistent with the fact that if $\gamma=(x,y)$ is a spacelike (resp.\ timelike) curve in $\lo^2$ such that $\kappa =\kappa (\rho)$, then
$\hat \gamma=(y,x)$ is a timelike (resp.\ spacelike) curve in $\lo^2$ such that $\kappa =\kappa (\rho)$.
Observe that we really only need to compute $\rho^+:=\rho^+_{1}=\rho^-_{-1}$ and $\rho^-:=\rho^+_{-1}=\rho^-_{1}$, and consequently only $\nu^+:=\nu^+_{1}=\nu^-_{-1}$ and $\nu^-:=\nu^+_{-1}=\nu^-_{+1}$.
}
\end{remark}
We show two illustrative examples applying steps (i)-(iii) in Remark~\ref{c rho}:
\begin{example}[$\kappa \equiv 0 $]
\label{lines r}
{\rm 
Then $\mathcal K\equiv c\in\R$, $s=\int \rho \, d\rho /\sqrt{c^2 \pm \rho^2}=\pm \sqrt{c^2\!\pm\!\rho^2} $.
So we get $\rho^+(s)=\sqrt{s^2-c^2} $, $|s| \geq c$, and  $\rho^-(s)=\sqrt{c^2-s^2} $, $|s| \leq c$. Accordingly, we have that $\nu^+ (s)=-\arccoth (s/c)  $ and $\nu^- (s)=\arctanh (s/c)  $. Then, a straightforward computation leads to the vertical geodesics $x=c$ at the spacelike case $\epsilon =1 $ and to the horizontal geodesics $y=c$ at the timelike case $\epsilon =-1 $. Up to the $\nu$-orthochrone Lorentz transformations $R_\nu$ given in \eqref{Rnu}, we reach in this way all the causal geodesics of $\lo^2$. 
When $\mathcal K\equiv 0$ we get exactly the geodesics $\alpha_{\varphi_0}$ and $\beta_{\phi_0}$ passing through the origin (see Figure \ref{Fig:geodesics}).}
\end{example}
\begin{example}[$\kappa\equiv 2k_0> 0$]
\label{circles r}
{\rm In this case, $\mathcal K (\rho)=k_0 \rho^2+c$, $c\in \R$, and 
$s=\int \rho \, d\rho /\sqrt{(k_0 \rho^2+c)^2\pm\rho^2} $, that in general is far from trivial. We consider the easiest case $\mathcal K(\rho)=k_0 \rho^2$, i.e.\ $c=0$ above. So we get $\rho^+(s)=\sinh (k_0 s)/k_0 $, $s\in \R$,  and  $\rho^-(s)=\cosh (k_0 s)/k_0 $, $s\in \R$.
In addition, in both cases $\nu^+(s)=\nu^-(s)=k_0 s$.  Using \eqref{pseudopolar}, we arrive (up to orthochrone transformations) to arbitrary pseudocircles of radius $1/2k_0$ passing through the origin (see Figure \ref{Fig:circlesO}). If $c\neq 0$, the integrals are much more complicated and the procedure is much longer and tedious to get the remaining pseudocircles in $\lo^2$.}
\end{example}
\begin{remark}\label{difficulties r}
{\rm The main difficulties one can find carrying on the strategy described in Theorem~\ref{quadratures rho}
(or in Remark~\ref{c rho}) to determine a Lorentzian curve whose curvature is $\kappa=\kappa(\rho)$ are the following:
\begin{enumerate}
\item The integration of $s=s(\rho)$:  Even in the case $\mathcal K(\rho)$ were polynomial,
the integral is not necessarily elementary. For example, when $\mathcal K(\rho)$ is a quadratic polynomial, it can be solved using Jacobian elliptic functions (see \cite{BF71}). This is equivalent to $\rho\,\kappa (\rho)$ be linear, so that $\kappa(\rho)=2\lambda + \mu/\rho$, $\lambda, \mu \in \R$. We will study such curves in Section \ref{Med spirals}.
\item The previous integration give us $s=s(\rho)$; it is not
always possible to obtain explicitly $\rho=\rho^+(s)$ and $\rho=\rho^-(s)$, what is necessary to determine the curve.
\item Even knowing explicitly $\rho=\rho^+(s)$ and $\rho=\rho^-(s)$, the integration to get $\nu^+ (s)$ and $\nu^- (s)$ may be impossible to perform using elementary or known functions.
\end{enumerate}
}
\end{remark}
Nevertheless, along the paper we will study different families where we are successful with the procedure described in Theorem~\ref{quadratures rho} and we will recover some known curves and find out new curves in $\lo^2$ characterized by this sort of geometric properties.

\subsection{Curves in $\lo^2$ such that $\kappa =\kappa(y\!-\!x)$}
\begin{figure}[h!]
\begin{center}
\includegraphics[height=4cm]{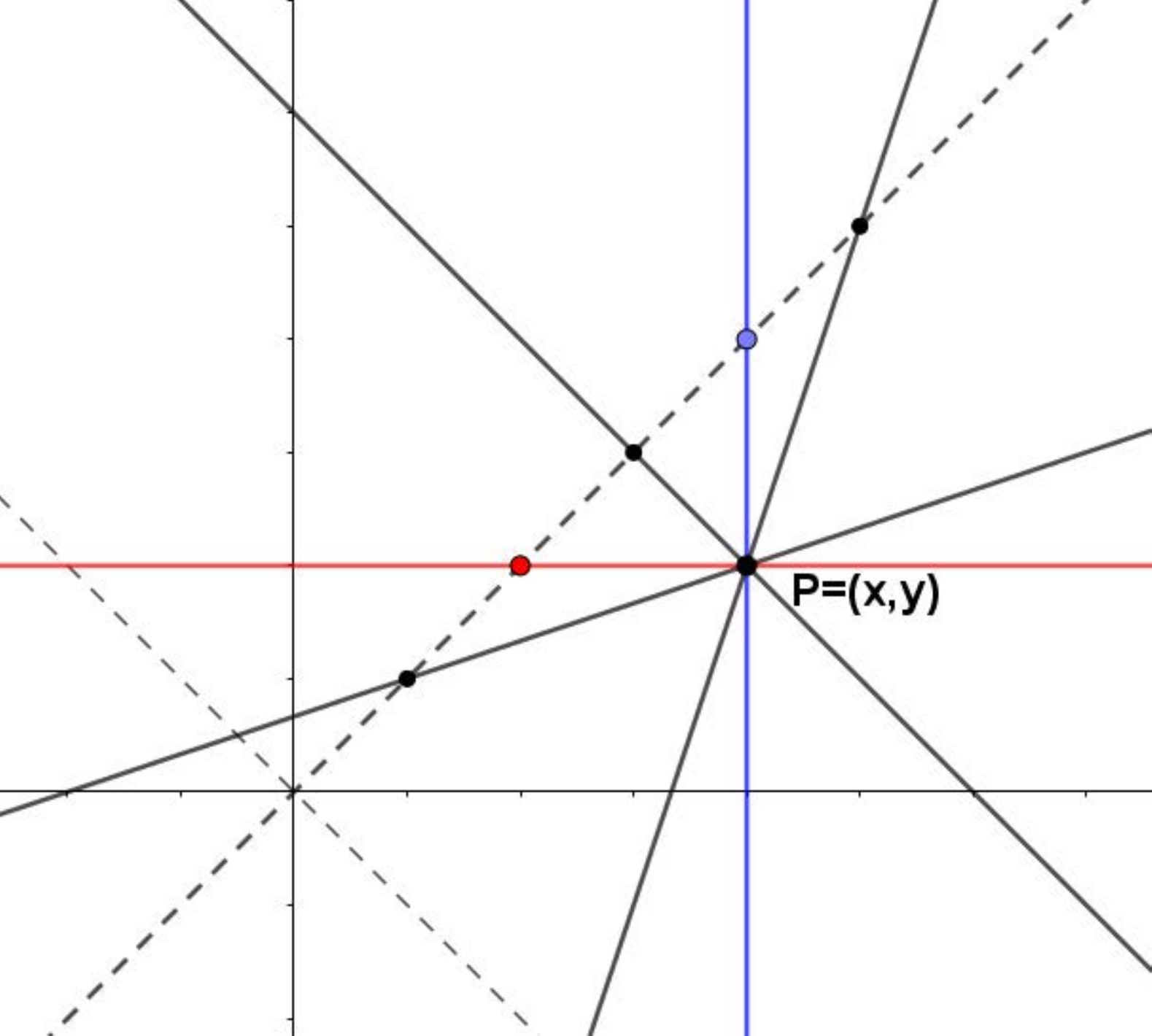}
\caption{Spacelike and timelike geodesics in $\lo^2$ passing through $P$.}
\label{Fig:inter k(v)}
\end{center}
\end{figure}

Given a spacelike or timelike curve $\gamma=(x,y)$ in $\lo^2$, we are now interested in the analytical condition $\kappa =\kappa (y-x)$. We look for its geometric interpretation. For this purpose,
we define the {\em Lorentzian pseudodistance between two points} by 
$$\delta: \lo^2 \times \lo^2 \rightarrow [0,+\infty), 
\ \delta (P,Q)=\sqrt{|g(\overrightarrow{PQ},\overrightarrow{PQ} )|}.$$ 
We fix the lightlike geodesic $x=y$. Given an arbitrary point $P\!=\!(x,y)\!\in\!\lo^2$, $x\neq y$, we consider all the spacelike and timelike geodesics $\gamma_m$ with slope $m\in \R \cup \{ \infty \}$, $m\neq 1$, passing through $P$, and let 
$P'=(\frac{mx-y}{m-1},\frac{mx-y}{m-1})$ the crossing point of $\gamma_m$ and the lightlike geodesic $x=y$ 
(see Figure \ref{Fig:inter k(v)}). 
Then:
$$ \delta(P,P')^2=(y-x)^2 \left| \frac{m+1}{m-1} \right|.$$
So $ \delta(P,P')^2=(y-x)^2$ if and only if $m=0$ or $m=\infty$.
Thus $|y-x|$ is the Lorentzian pseudodistance from $P\!=\!(x,y)\in\lo^2$, $x\neq y$, to the lightlike geodesic $x=y$
through the horizontal timelike geodesic or the vertical spacelike geodesic.

For the sake of simplicity, we introduce new coordinates $(u,v)$, $u,v\in \R$, in $\lo^2$ by the equations
\begin{equation}\label{uv}
u=y+x, \ v=y-x \Leftrightarrow x=\frac{u-v}{2}, \ y = \frac{u+v}{2}.
\end{equation}
Then we are considering the condition $\kappa = \kappa (v)$.

At a given point $\gamma (s)=(x(s),y(s))$ on the curve, the {\em geometric linear momentum (with respect to the lightlike geodesic $x=y$ or $u$-axis)} $\mathcal K$  is given by
\begin{equation}\label{u-momentum}
\mathcal K(s) = \dot u(s).
\end{equation}
In physical terms, using Noether's Theorem, $\mathcal K$ may be interpreted as the linear momentum with respect to the $u$-axis of a particle of unit mass with unit-speed and trajectory $\gamma $. 

Since  $\gamma$ is unit-speed, we have simply that $ \dot u \, \dot v= \epsilon$. Using \eqref{u-momentum}, 
we easily obtain that
\begin{equation}\label{differentials v}
ds = \epsilon \mathcal K dv, \quad du = \mathcal K ds.
\end{equation}
Thus, given $\mathcal K=\mathcal K (v)$ as an explicit function, looking at \eqref{differentials v} one may attempt to compute $v(s)$ and $u (s) $ in three steps: integrate to get $s=s(v)$, invert to get $v=v(s)$ and integrate to get $u =u(s)$. 

In addition, we have that the curvature $\kappa $ satisfies \eqref{FrenetEq1}, that implies $\ddot u = \kappa \dot u $ and so $ \dot {\mathcal K} = \kappa\, \mathcal K$ by taking into account \eqref{u-momentum}. From \eqref{differentials v} we deduce that $d\mathcal K=\epsilon \kappa \, \mathcal K^2 \, dv$ and, since we are assuming that $\kappa =\kappa (v)$, we finally arrive at 
\begin{equation}\label{K primitive v}
-\epsilon \, d (1/\mathcal K)  =  \kappa(v)dv,
\end{equation}
that is, $-\epsilon / \mathcal K (v)$ can be interpreted as an anti-derivative of $\kappa (v)$.

As a summary, we have proved the following result in the spirit of Theorem 3.1 in \cite{S99}.
\begin{theorem}\label{quadratures v}
Let $\kappa=\kappa(v)$ be a continuous function. 
Then the problem of determining locally a spacelike or timelike curve in $\lo^2$ whose curvature is $\kappa(v)$ with geometric linear momentum $\mathcal K (v)$ satisfying \eqref{K primitive v} ---$|v|$ being the (non constant) Lorentzian pseudodistance through the horizontal or vertical geodesics to the $u$-axis--- is solvable by quadratures considering the unit speed curve $(u(s),v(s))$, 
where $v(s)$ and $u (s)$ are obtained through (\ref{differentials v}) after inverting $s=s(v)$,
being $\epsilon=1$ at the spacelike case and $\epsilon=-1$ at the timelike case.
Such a curve is uniquely determined by $\mathcal K (v)$ up to a translation in the $u$-direction (and a translation of the arc parameter $s$).
\end{theorem}
\begin{remark}\label{c v}
{\rm If we prescribe $\kappa=\kappa(v)$, the method described in Theorem~\ref{quadratures v} clearly implies the computation of three quadratures, following the sequence:
\begin{enumerate}
\item[\rm (i)]  Anti-derivative of $\kappa (v)$: $$\int \! \kappa(v)dv= -\epsilon / \mathcal K (v). $$
\item[\rm (ii)]  Arc-length parameter $s$ of $(u(s),v(s))$ in terms of $v$: $$s=s(v)=\epsilon \int \mathcal K (v)\,dv,$$
and inverting $s=s(v)$ to get $v=v(s)$.
\item[\rm (iii)] First coordinate of $(u(s),v(s))$ in terms of $s$: $$u(s)=\int  K (v(s))ds.$$
\end{enumerate}
We note that we get a one-parameter family of curves in $\lo^2$ satisfying $\kappa=\kappa(v)$ according to the geometric linear momentum chosen in (i). It will distinguish geometrically the curves inside a same family by their relative position with respect to the $u$-axis. We remark that we can recover $\kappa$ from $\mathcal K$ by means of $\kappa (v)=  -\epsilon d(1/\mathcal K)/dv$. In addition, we remark that if $(u(s),v(s))$ is the curve corresponding to $\epsilon =1$, then $(-u(s),v(s))$ is the curve corresponding to $\epsilon =-1$; they have geometric linear momentum with opposite sign, but the same intrinsic equation $\kappa =\kappa (s)$.
}
\end{remark}
We show two illustrative examples applying steps (i)--(iii) in Remark~\ref{c v}:
\begin{example}[$\kappa \equiv 0 $]
\label{lines v}
{\rm 
Then $\mathcal K\equiv -\epsilon /c$, $c\neq 0$. It is very easy to get $v(s)=-c \,s $ and $u(s)=-\epsilon\,s/c$, $s \in \R$, that parametrize the line passing through the origin with slope $m=\frac{\epsilon + c^2}{\epsilon-c^2}$. We observe that $c=0$ implies $m=1$, corresponding to the lightlike geodesic given by $u$-axis. If $\epsilon =1$, then $|m|>1$ and we obtain the spacelike geodesics $\alpha_{\varphi_0}$. If $\epsilon =-1$, then $|m|<1$ and we have the timelike geodesics $\beta_{\phi_0}$. See Figure \ref{Fig:geodesics}.}
\end{example}
\begin{example}[$\kappa\!\equiv\! k_0\!> \!0$]
\label{circles v}
{\rm  Now $\mathcal K (v)=-\epsilon/(c+ k_0 v)$, $c\in\R$. Then it is not difficult to get that $v(s)=(e^{-k_0 s}-c)/k_0$ and $u(s)=-\epsilon e^{k_0 s}/k_0$. Using \eqref{uv}, 
if $\epsilon =1$, we get $ x(s)=(-\cosh(k_0 s)+c/2)/k_0$ and $ y(s)=-(\sinh(k_0s)+c/2)/k_0 $;
and if $\epsilon =-1$, we get $ x(s)=(\sinh(k_0 s)+c/2)/k_0$ and $ y(s)=(\cosh(k_0s)-c/2)/k_0 $.
They correspond respectively to spacelike and timelike pseudocircles in $\lo^2$ of radius $1/k_0$ (see Figure \ref{Fig:circlesbis}).
When $c=0$, we obtain the branches of $x^2-y^2=\epsilon /k_0^2$ (travelled with positive curvature $k_0$), that are asymptotic to the light cone of $\lo^2$.}
\begin{figure}[h!]
\begin{center}
\includegraphics[height=4cm]{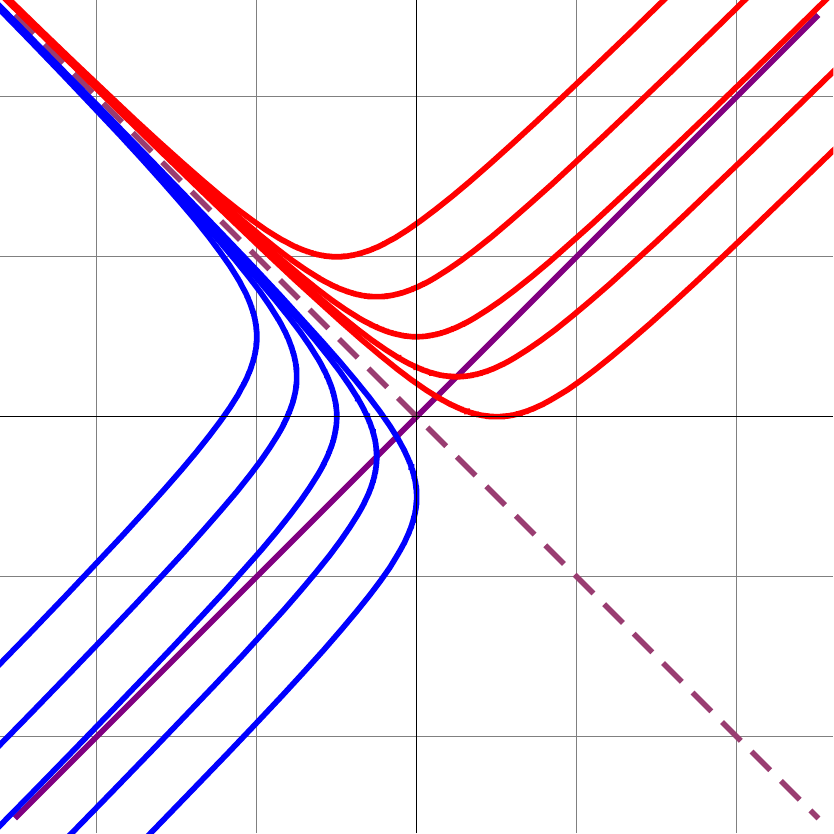}
\caption{Spacelike (blue) and timelike (red) pseudocircles in $\lo^2$ with constant positive curvature.}
\label{Fig:circlesbis}
\end{center}
\end{figure} 
\end{example}

\vspace{0.3cm}

\section{The Lorentzian Norwich spiral and its generalizations}\label{Sturm}
\begin{figure}[h!]
\begin{center}
\includegraphics[height=4cm]{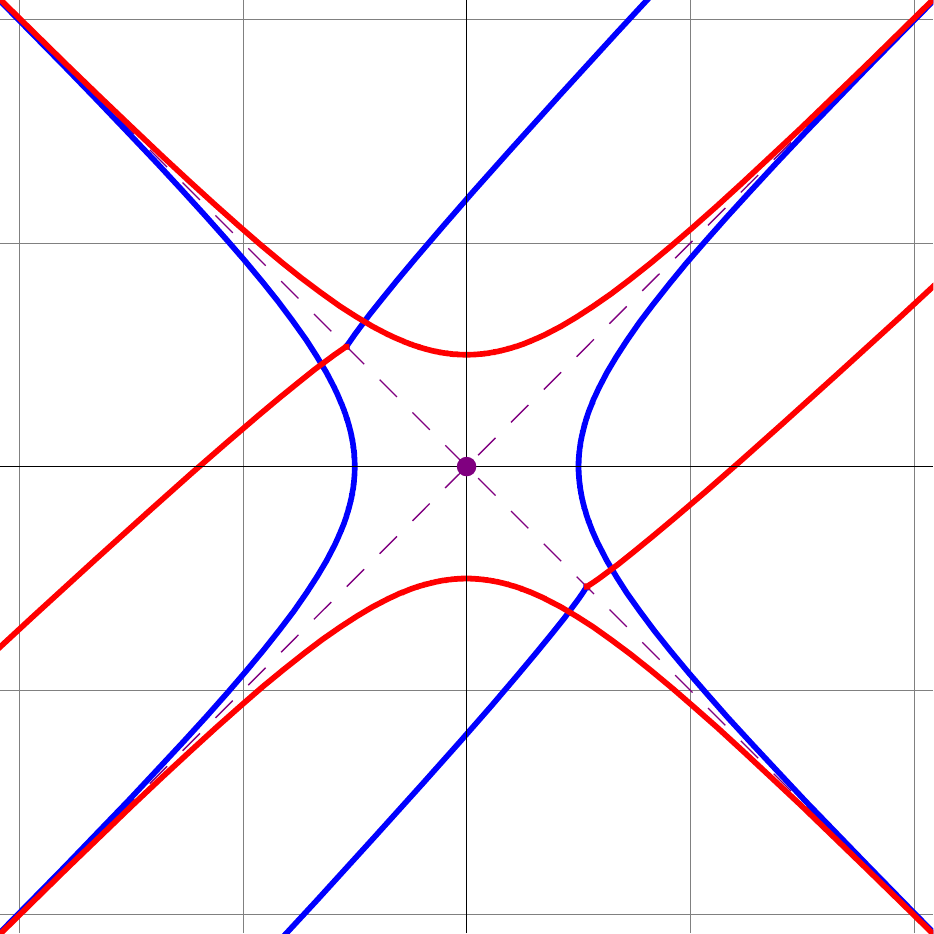}
\caption{Lorentzian Norwich spiral.}
\label{Fig:Norwich}
\end{center}
\end{figure}

In the context of Euclidean Geometry, the Norwich spiral is the curve (other than a circle) such that the radius of curvature $\mathcal R := 1/\kappa$ at any point is equal to the distance to a fixed point.  It is also known as the Sturm spiral. We aim to find out the Lorentzian version of this interesting curve, using the notion of Lorentzian pseudodistance $\rho$. 
Up to a translation, we can take the origin as the fixed point and hence we are considering the geometric condition
\begin{equation}\label{Norwich}
\kappa (\rho)=\frac{1}{\rho} \Leftrightarrow \mathcal R (\rho)=\rho, \ \rho >0,
\end{equation}
which is invariant under dilations. 
Of course we also have the trivial solutions of \eqref{Norwich} corresponding to the pseudocircles $\rho \equiv \rho_0 >0$.
For the non-trivial ones, we follow the steps described in Remark \ref{c rho}. Using \eqref{Norwich}, we have:
$$\int \! \rho\, \kappa(\rho)d\rho= \rho+c, \ c\neq 0. $$
In addition: 
$$s=s(\rho)=\int\frac{\rho \, d\rho}{\sqrt{(\rho +c)^2\pm \rho^2}},$$
which implies, using plus sign, that
$$  s = \frac{\sqrt{2(\rho^+)^2+2c\rho^++c^2}}{2}-\frac{c}{2\sqrt 2} \arcsinh (2\rho^+/c+1)$$
and, using minus sign, that
$$ s=\frac{\rho^--c}{3c}\sqrt{2c\rho^- +c^2}.$$
In both cases, it is not possible to invert $s=s(\rho)$ to get $\rho=\rho^+(s)$ and $\rho=\rho^-(s)$ in an explicit way (see Remark \ref{difficulties r}). But if we use a new parameter coming from $ds=\rho\,dt$, we easily arrive at 
$$ \rho^+(t)=\frac{c}{2}\left( \sinh (\sqrt 2 t)-1 \right), \, t >\frac{1}{\sqrt 2} \arcsinh 1,$$
and
$$ \rho^-(t)=\frac{c}{2}(1-t^2), \, |t|<1,$$
respectively. The relation of $t$ with the arc-length parameter $s$ is given by 
$$ s=\frac{c}{2} \left(  \frac{\cosh (\sqrt 2 t)}{\sqrt 2}-t \right), \quad s=\frac{c}{2} \left(t-\frac{t^3}{3}\right), $$
respectively.
Using again the parameter $t$, we have that $d\nu^\pm = (1+c/\rho^\pm)dt$ and we finally get that
$$ \nu^+ (t)=t+\log \left(  \frac{\sinh (\frac{\sqrt 2 \, t -\arcsinh 1}{2}) }{\cosh (\frac{\sqrt 2 \, t +\arcsinh 1}{2}) } \right),
\ t >\frac{1}{\sqrt 2} \arcsinh 1, $$
and
$$\nu^-(t)=t+2 \arctanh t, \ |t|<1.  $$
Using the above expressions in \eqref{pseudopolar}, we obtain the explicit parametrizations for a spacelike and timelike curve in $\lo^2$ satisfying \eqref{Norwich}.
It is quite remarkable that different values of the constant $c$ only produce homothetic curves in this case; recall that \eqref{Norwich} is invariant under dilations. We call this curve the Lorentzian Norwich spiral (see Figure \ref{Fig:Norwich}).

In \cite{IUM15}, the authors studied the Lorentzian curves satisfying the geometric condition
\begin{equation}\label{Sturmian}
\kappa (\rho)=\frac{\mu}{\rho}, \ \mu >0.
\end{equation}
The Norwich (or Sturm) spiral obviously corresponds to the case $\mu =1$.
Making use of the co-moving frame associated to a curve in $\lo^2$, the cases $\mu >1$  and $0<\mu <1$ (using the parameter $t$ given by $ds=\rho dt$)  were completely determined in \cite{IUM15} in terms of elementary functions. 

\vspace{0.3cm}

\section{Curves in $\lo^2$ such that $\kappa (\rho)=2\lambda + \mu/\rho$}\label{Med spirals}
The aim of this section is to study the curves in the Lorentz-Minkowski plane satisfying the geometric condition
$\kappa (\rho)=2\lambda + \mu/\rho$,
where $\rho>0$ denotes the Lorentzian pseudodistance from the origin and $\lambda, \mu \in \R$. 
If $\mu=0$, then $\kappa$ is constant and we get only geodesics and pseudocircles. 
If $\lambda =0$, then $\kappa = \mu /\rho$ and we get the generalized Sturm spirals studied in
\cite{IUM15} (see Section \ref{Sturm}). Thus we afford only the case $\lambda > 0, \mu \neq 0$. 
Our purpose will be to find analytic representations for these curves using elementary functions, such as it happened to the generalized Sturmian spirals.
Up to dilations, we can also consider $\lambda = 1$ without loss of generality. 
In fact, if $\gamma$ is a Lorentzian curve
satisfying $\kappa_\gamma =2\lambda + \mu / |\gamma|$, $\lambda >0, \mu\neq 0$, 
then $\kappa_{a\gamma} =2\lambda
/a + \mu / |a \gamma|$, for any $a >0$. 
Taking $a=\lambda$ it is enough to consider the condition
\begin{equation}\label{condition spi}
\kappa (\rho)=2+\frac{\mu}{\rho}, \, \mu \neq 0,
\end{equation}
since if $\gamma$ satisfies (\ref{condition spi}) then $ \gamma /\lambda$ will verify $\kappa (\rho) = 2\lambda + \mu/\rho$.

If we look for trivial solutions of (\ref{condition spi}), that is, curves in $\lo^2$ with constant curvature, we arrive at the pseudocircles $\rho\equiv (1-\mu)/2$, $\forall\mu <1$. In the non-trivial case, following Remark \ref{c rho} and using \eqref{condition spi}, we have:
\begin{equation}\label{signo}
\int \! \rho\,\kappa(\rho)d\rho= \rho^2+\mu \rho+c. 
\end{equation}
Hence we must consider the geometric angular momentum $\mathcal K(\rho)=\rho^2 +\mu \rho+ c$, $c\in \R$. Then:
\begin{equation}\label{r_general_1}
s = \int \! \frac{\rho \, d\rho}{\sqrt{( \rho^2 + \mu \rho +c )^2 \pm \rho^2}}
\end{equation}
and
\begin{equation}\label{arg_general}
\nu(s)= s + \mu \int \frac{ds}{\rho(s)}+ c \int \frac{ds}{\rho(s)^2} .
\end{equation}
We point out that if $c\neq 0$, \eqref{r_general_1} can be solved using Jacobi elliptic functions for $\rho=\rho(s)$ but \eqref{arg_general} involves very complicated elliptic integrals of second and third kind for $\nu = \nu (s)$ depending on the values of $c$. Therefore we only pay attention to the nicest case $c=0$, corresponding (see \eqref{signo}) to $$\mathcal K(\rho)=\rho^2 +\mu \rho.$$
Clearly (\ref{r_general_1}) reduces to
$$ s= \int \frac{d\rho}{\sqrt{( \rho+\mu)^2\pm 1}} $$
which gives, according to plus or minus sign,
\begin{equation}\label{r_explicit +}
\rho^+_\mu(s)=\sinh s - \mu , \ s > \arcsinh \mu,
\end{equation}
and
\begin{equation}\label{r_explicit -}
\rho^-_\mu(s)=\cosh s - \mu, \ \left\{
\begin{array}{l}
\ s\in \R, {\rm \ when \ } \mu < 1 \\ |s| > \arccosh \mu, {\rm \ when \ } \mu \geq 1
\end{array}
\right.
\end{equation}
respectively.
Accordingly, (\ref{arg_general}) translates into
\begin{equation}\label{theta_explicit +}
\nu^+_\mu(s)= s + \mu \int \frac{ds}{\sinh s -\mu }.
\end{equation}
and
\begin{equation}\label{theta_explicit -}
\nu^-_\mu(s)= s + \mu \int \frac{ds}{\cosh s -\mu }.
\end{equation}
For the integration of \eqref{theta_explicit +}, we write $\mu=\sinh \eta$, $\eta \in \R$, and so $\rho^+_\eta(s)=\sinh s - \sinh \eta $, $ s > \eta$, and then: 
\begin{equation}\label{nu +}
\nu^+_\eta(s)=s+ \tanh \eta \, \log  \left(  \frac{\sinh (\frac{s-\eta}{2})}{\cosh (\frac{s+\eta}{2})} \right), \ s> \eta.
\end{equation}
We must distinguish cases according to the values of $\mu$ to perform the integration of (\ref{theta_explicit -}), in order to
use its result joint to \eqref{nu +} in \eqref{pseudopolar}, obtaining in this way the explicit unit-speed parametrizations of a spacelike and timelike curve in $\lo^2$ satisfying \eqref{condition spi} with geometric angular momentum $\mathcal K(\rho)=\rho^2 +\mu \rho$.
\subsection{Case $\mu = 1$} 
\begin{figure}[h!]
\begin{center}
\includegraphics[height=4cm]{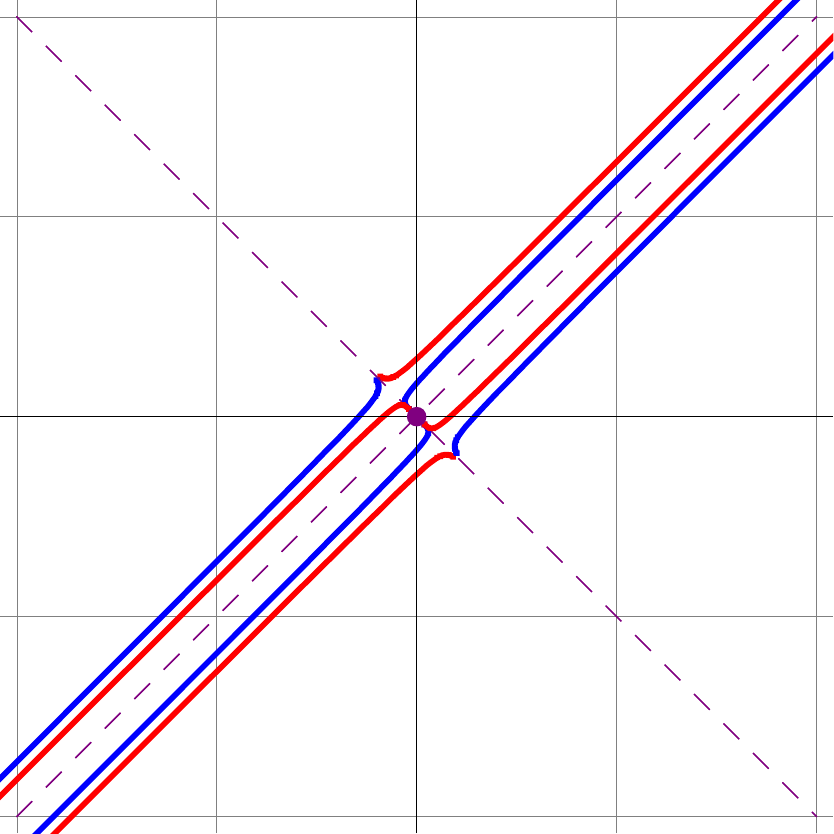}
\caption{Spacelike (blue) and timelike (red) curve such that $\mathcal K(\rho)=\rho^2 + \rho$.}
\label{Fig:muplus1}
\end{center}
\end{figure}

Then we have $\rho^-_1(s)=\cosh s - 1 $, $s\neq 0$, and 
\begin{equation}\label{nu- mu1}
\nu^-_1(s)=s-\coth (s/2), \, s\neq 0.
\end{equation}
See Figure \ref{Fig:muplus1}.

\subsection{Case $\mu = - 1$} 
\begin{figure}[h!]
\begin{center}
\includegraphics[height=4cm]{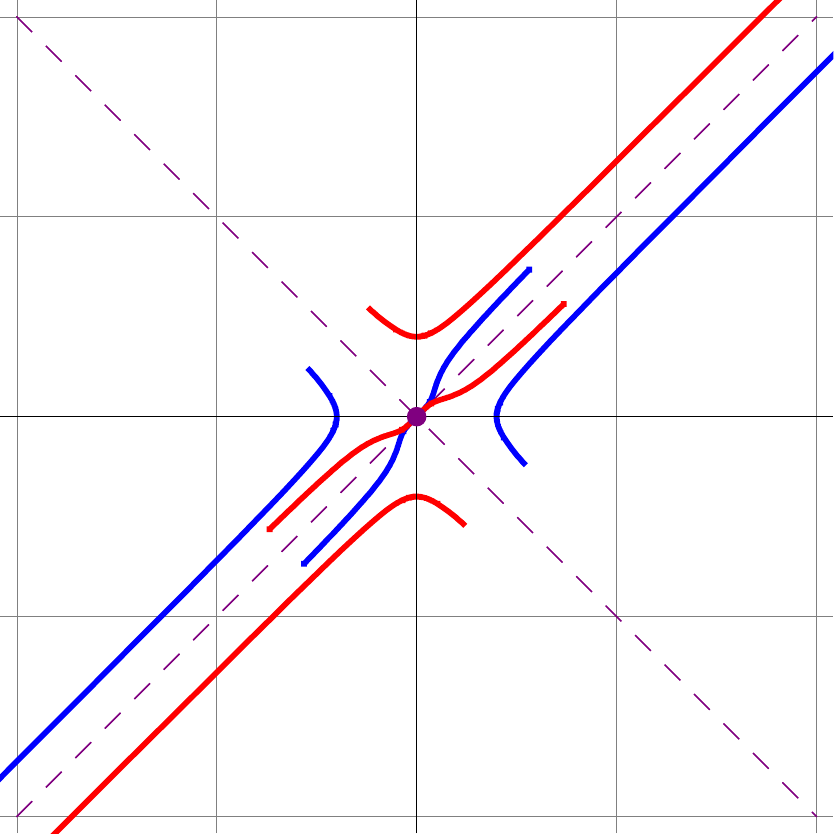}
\caption{Spacelike (blue) and timelike (red) curve such that $\mathcal K(\rho)=\rho^2 - \rho$.}
\label{Fig:muminus1}
\end{center}
\end{figure}

Now $\rho^-_{-1}(s)=\cosh s +1 $, $s \in \R $, and
\begin{equation}\label{nu- mu-1}
\nu^-_{-1}(s)=s-\tanh (s/2), \, s\in \R.
\end{equation}
See Figure \ref{Fig:muminus1}.

\subsection{Case $|\mu | < 1$}  
\begin{figure}[h!]
\begin{center}
\includegraphics[height=4cm]{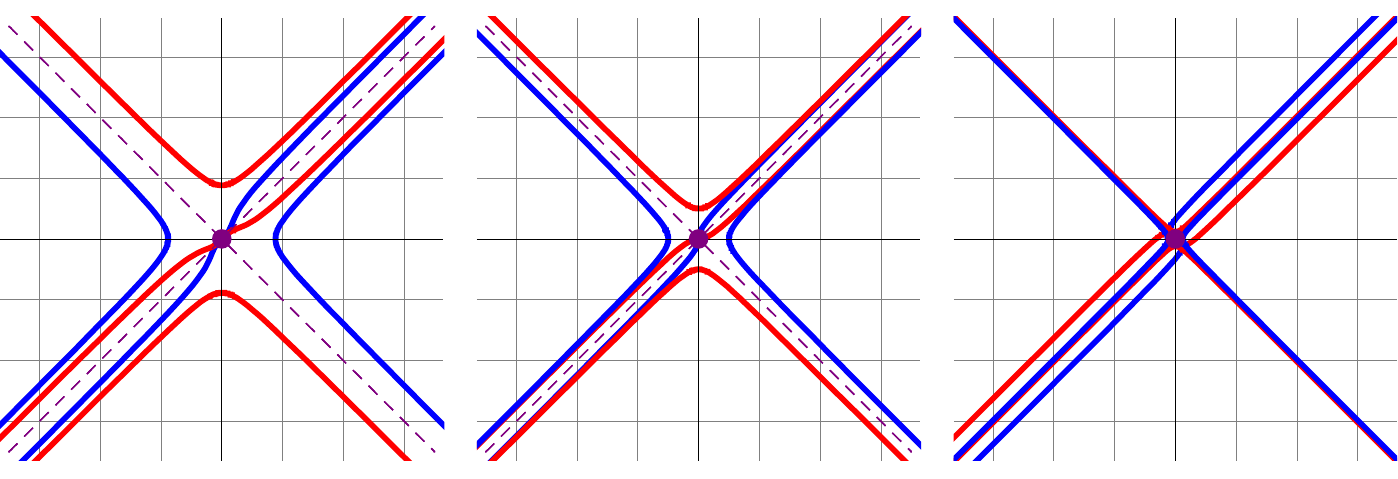}
\caption{Spacelike (blue) and timelike (red) curves such that $\mathcal K(\rho)=\rho^2 + \cos \alpha \,  \rho$, $\alpha=\frac{\pi}{4}, \frac{\pi}{2} (pseudocircle), \frac{3\pi}{4} $.}
\label{Fig:mualpha}
\end{center}
\end{figure}

In this case, we put $\mu =\cos \alpha$, with $0<\alpha <\pi$. 
Then  $\rho^-_\alpha(s)=\cosh s - \cos \alpha $, $s\in \R$, and
\begin{equation}\label{nu- mu<1}
\nu^-_\alpha(s)=s+2\cot \alpha \, \arctan \left(  \cot (\alpha/2)\tanh (s/2) \right), \, s\in \R. 
\end{equation}
See Figure \ref{Fig:mualpha}.

\subsection{Case $\mu  > 1$} 
\begin{figure}[h!]
\begin{center}
\includegraphics[height=4cm]{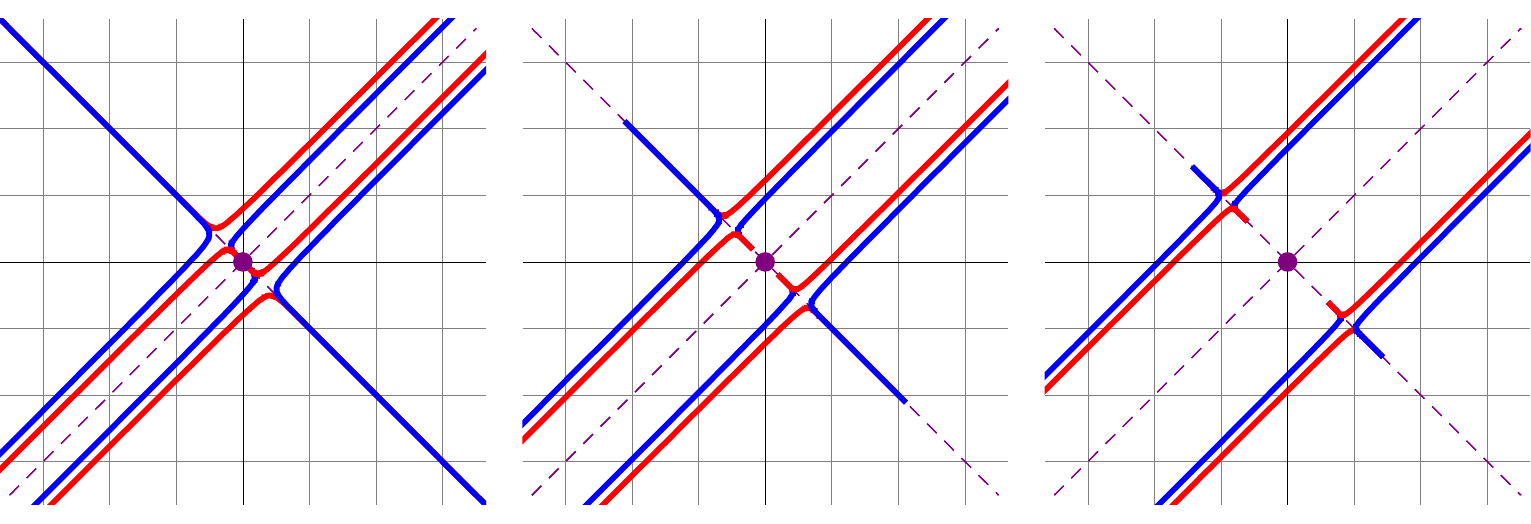}
\caption{Spacelike (blue) and timelike (red) curves such that $\mathcal K(\rho)=\rho^2 +  \cosh \delta \, \rho$, $\delta >0$.}
\label{Fig:mudelta}
\end{center}
\end{figure}

Now we write $\mu=\cosh \delta$,  $\delta >0$. Then $\rho^-_\delta(s)=\cosh s - \cosh \delta $, $|s|>\delta$, and 
\begin{equation}\label{nu- mu>1}
\nu^-_\delta(s)=s+ \coth \delta \, \log  \left(  \frac{\sinh (\frac{s-\delta}{2})}{\sinh (\frac{s+\delta}{2})} \right),  \,  |s|> \delta. 
\end{equation}
See Figure \ref{Fig:mudelta}.

\subsection{Case $\mu  < -1$} 
\begin{figure}[h!]
\begin{center}
\includegraphics[height=4cm]{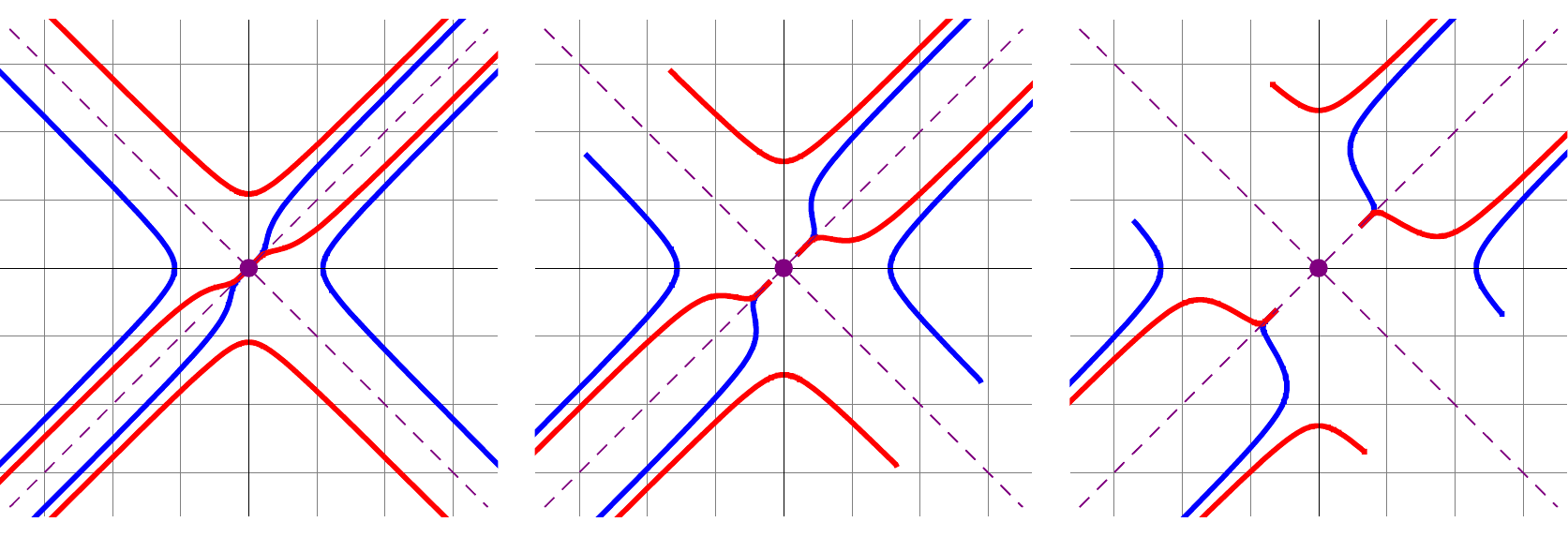}
\caption{Spacelike (blue) and timelike (red) curves such that $\mathcal K(\rho)=\rho^2 - \cosh \tau \,  \rho$, $\tau >0$.}
\label{Fig:mutau}
\end{center}
\end{figure}

If we put $\mu =-\cosh \tau$, $\tau >0$, we have that $\rho^-_\tau(s)=\cosh s + \cosh \tau $, $s\in \R$, and
\begin{equation}\label{nu- mu<-1}
\nu^-_\tau(s)=s+ \coth \tau \, \log  \left(  \frac{\cosh (\frac{s-\tau}{2})}{\cosh (\frac{s+\tau}{2})} \right), \, s\in \R. 
\end{equation}
See Figure \ref{Fig:mutau}.

\vspace{0.3cm}

\section{Lorentzian sinusoidal spirals}
In this section, we want to study those spacelike and timelike curves in $\lo^2$ such that 
$$\kappa (\rho)=\lambda \, \rho^{n-1},\, \lambda >0, \,n\in \R.$$
If $n=1$, we get only pseudocircles.  If $n=0$, we arrive at the generalized Sturmian spirals (cf.\ \cite{IUM15}).
So we start considering $n\neq 0$. We follow the strategy described in Remark \ref{c rho} considering,
when $n\neq -1$, the geometric angular momentum given by $$\mathcal K(\rho)=  \frac{\lambda}{n+1}\rho^{n+1}.$$
If we combine steps (ii) and (iii) in Remark~\ref{c rho}, eliminating $ds$, we arrive at
$$
\frac{d\nu_+}{d\rho_+}
=\frac{\frac{\lambda}{n+1}\, \rho^{n-1}}{\sqrt{\frac{\lambda^2}{(n+1)^2} \rho^{2n}+1}},
\quad
\frac{d\nu_-}{d\rho_-}
=\frac{\frac{\lambda}{n+1}\, \rho^{n-1}}{\sqrt{\frac{\lambda^2}{(n+1)^2} \rho^{2n}-1}}.
$$
In this way, we deduce the pseudopolar equations of this family given by 
\begin{equation}\label{LorSinSpirals}
\begin{array}{c}
 \lambda \, \rho_+^{n} = (n+1)\sinh \left( n \nu_+ \right), \, n\neq 0, \, n\neq -1, \\ \\
 \lambda \, \rho_-^{n} = (n+1)\cosh \left( n \nu_- \right), \, n\neq 0, \, n\neq -1,
\end{array}
\end{equation}
that combined with \eqref{pseudopolar}, taking into account Remark \ref{c rho}, provide us the different curves of the family.
Taking into account Section 7.1  in \cite{CCIs17}, we will refer these curves as {\em Lorentzian sinusoidal spirals}. In this way, we deduce the following characterization of this wide family of curves in $\lo^2$.
\begin{corollary}\label{cor:sinusoidal}
The Lorentzian sinusoidal spiral \eqref{LorSinSpirals} is the only curve (up to $\nu$-orthochrone transformations) in $\lo^2$ with geometric angular momentum $\mathcal K(\rho)=  \frac{\lambda}{n+1}\rho^{n+1} $ (and curvature  $\kappa (\rho)=\lambda \, \rho^{n-1},\, \lambda >0$), $ n\neq 0, \, n\neq -1$.
\end{corollary}

\begin{figure}[h!]
\begin{center}
\includegraphics[height=3cm]{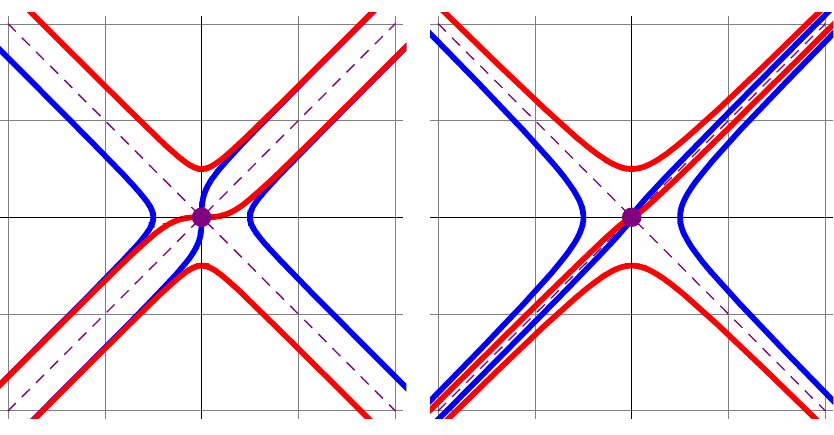}
\caption{Lorentzian sinusoidal spirals with $n=2$ (left) and $n=1/2$ (right).}
\label{Fig:LemnisCard}
\end{center}
\end{figure}

\begin{figure}[h!]
\begin{center}
\includegraphics[height=3cm]{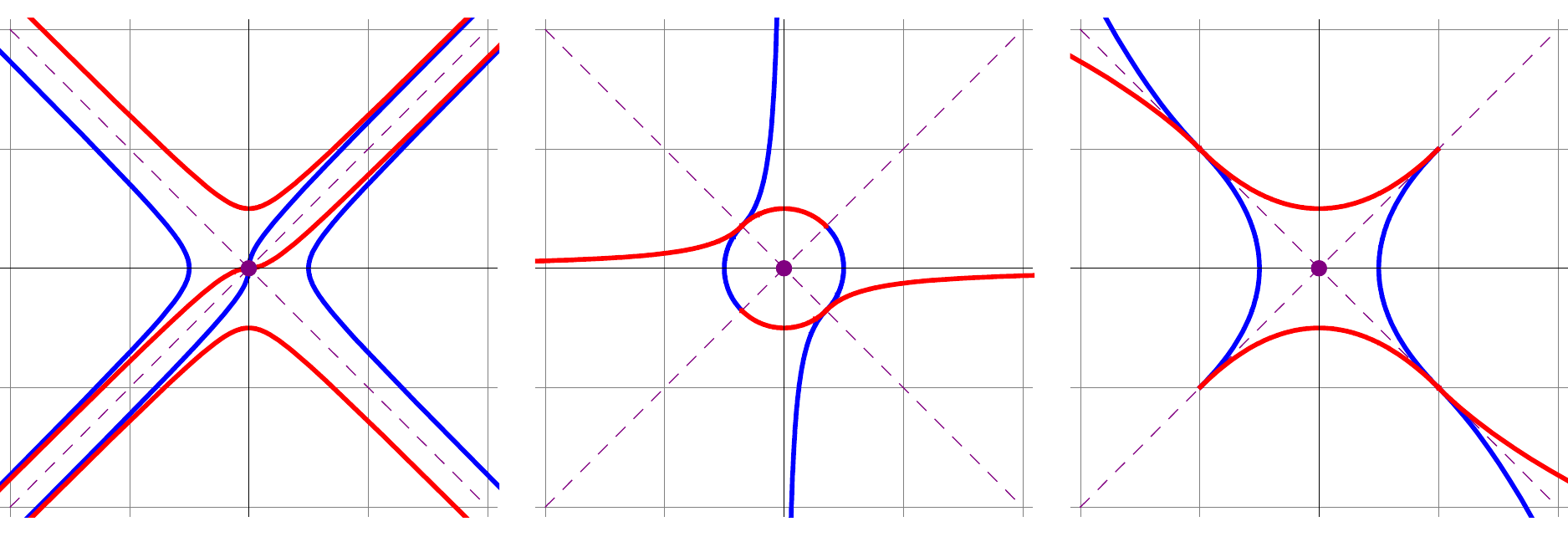}
\caption{Lorentzian sinusoidal spirals with $n=1$ (left), $n=-2$ (center) and $n=-1/2$ (right).}
\label{Fig:CirHipPar}
\end{center}
\end{figure}

\begin{figure}[h!]
\begin{center}
\includegraphics[height=4cm]{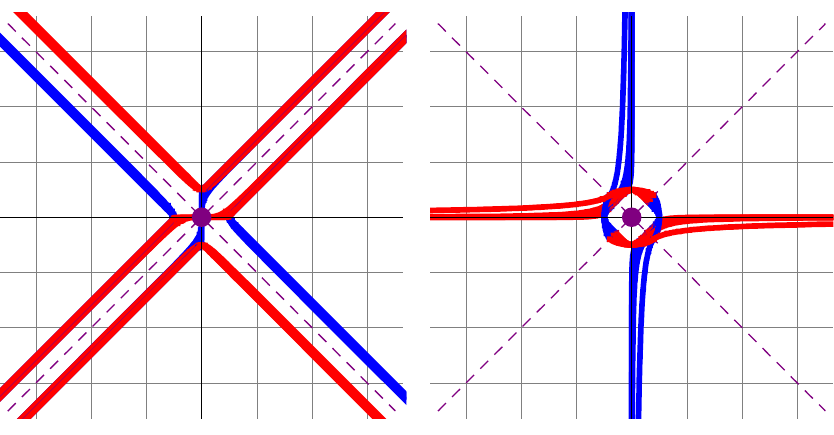}
\caption{Lorentzian sinusoidal spirals with $n \geq 5/2$ (left) and $n\leq -3/2$ (right), $n\in \Q$.}
\label{Fig:nPosNeg}
\end{center}
\end{figure}

\begin{remark}\label{re:conics}
{\rm Taking into account what happens in the Euclidean case (see Remark 7.1 in \cite{CCIs17}), the Lorentzian sinusoidal spirals take in the Lorentzian versions of some very interesting plane curves, including some conics. In particular, up to dilations, we emphasize the following curves:
\begin{enumerate}
\item[(i)] $n=2$: the {\em Lorentzian Bernoulli pseudolemniscate} defined by $$\rho_+^2= \sinh 2 \nu_+, \, \rho_-^2= \cosh 2 \nu_- $$ 
with $\mathcal K (\rho)= \rho^3$ (see Figure \ref{Fig:LemnisCard});
\item[(ii)] $n=1/2$: the {\em Lorentzian pseudocardioid} defined by $$\sqrt \rho_+= \sinh (\nu_+/2), \, \sqrt \rho_-=\cosh (\nu_-/2) $$
with $\mathcal K (\rho)= \rho^{3/2}$ (see Figure \ref{Fig:LemnisCard});
\item[(iii)] $n=1$: the pseudocircles given by $$\rho_+= \sinh  \nu_+, \, \rho_-= \cosh  \nu_- $$ 
with $\mathcal K (\rho)= \rho^2$ (see Figure \ref{Fig:CirHipPar});
\item[(iv)] $n=-2$: the {\em Lorentzian equilateral pseudohyperbolas} defined by $$\rho_+^2= -1/\sinh 2 \nu_+, \, \rho_-^2= 1/ \cosh 2 \nu_- $$
with $\mathcal K (\rho)= 1/\rho$ (see Figure \ref{Fig:CirHipPar});
\item[(v)] $n=-1/2$: the {\em Lorentzian pseudoparabolas} defined by $$\sqrt \rho_+=-1/ \sinh (\nu_+/2), \, \sqrt \rho_-=1/ \cosh (\nu_-/2) $$
with $\mathcal K (\rho)= \sqrt \rho$ (see Figure \ref{Fig:CirHipPar}).
\end{enumerate}
}
\end{remark}
If $n $ is rational, then the Lorentzian sinusoidal spirals could be considered Lorentzian versions of algebraic curves (see Figure \ref{Fig:nPosNeg}).

\vspace{0.3cm}

\section{Curves in $\lo^2$ such that $\kappa (v)=a v + b$, $a \neq 0, \, b \in \R $}
In this section we will study those spacelike and timelike curves in $\lo^2$ satisfying
\begin{equation}\label{cond_genelasticae}
\kappa(v)  = a v + b, \ a \neq 0, \, b \in \R,
\end{equation}
and we will show its close relationship with a certain class of elastic curves of $\lo^2$.
Recall that a unit speed spacelike or timelike curve $\gamma$ in $\lo^2$ is said to be an {\em elastica under tension $\sigma$} (see \cite{SAEY11}) if it satisfies the differential equation
\begin{equation}\label{eq_elastica}
2 \ddot \kappa - \kappa^3 -\sigma \kappa =0,
\end{equation}
for some value of $\sigma \in \R$.
Multiplying
(\ref{eq_elastica}) by $2\dot \kappa$ and integration allow us to introduce the {\em energy} $E\in\R$ of an elastica:
\begin{equation}\label{energy_elastica}
E:=  \dot \kappa ^2 - \frac{1}{4}\kappa^4 -\frac{\sigma}{2} \kappa^2.
\end{equation}
Given $\gamma =(u,v)$ satisfying \eqref{cond_genelasticae} with $a>0$ without restriction, we take $\hat \gamma = \sqrt{
a/2} (u,v+ b/a)$ and then, up to a translation in the $v$-direction and a dilation, we can only afford the condition
\begin{equation}\label{cond_elastica}
\kappa (v) = 2 v.
\end{equation}
We follow the strategy described in Remark~\ref{c v} in order to control the spacelike or timelike curves 
$(u(s),v(s))$ in $\lo^2$ satisfying (\ref{cond_elastica}). 
First, we obtain that the geometric linear momentum is given by
\begin{equation}\label{Kelas}
\mathcal K(v)=-\frac{\epsilon }{v^2+c}, \,  c\in \R.
\end{equation}
The arc-parameter can be deduced from
\begin{equation}\label{s elas}
s=s(v)=- \int \frac{dv}{v^2+c},
\end{equation}
and after inverting to get $v=v(s)$ and using \eqref{Kelas}, we have that
\begin{equation}\label{u elas}
u(s)=-\epsilon\int \frac{ds}{v(s)^2+c}.
\end{equation}
Then we must distinguish cases according to the values of $c\in \R$.
\subsection{Case $c=0$}
\begin{figure}[h!]
\begin{center}
\includegraphics[height=4cm]{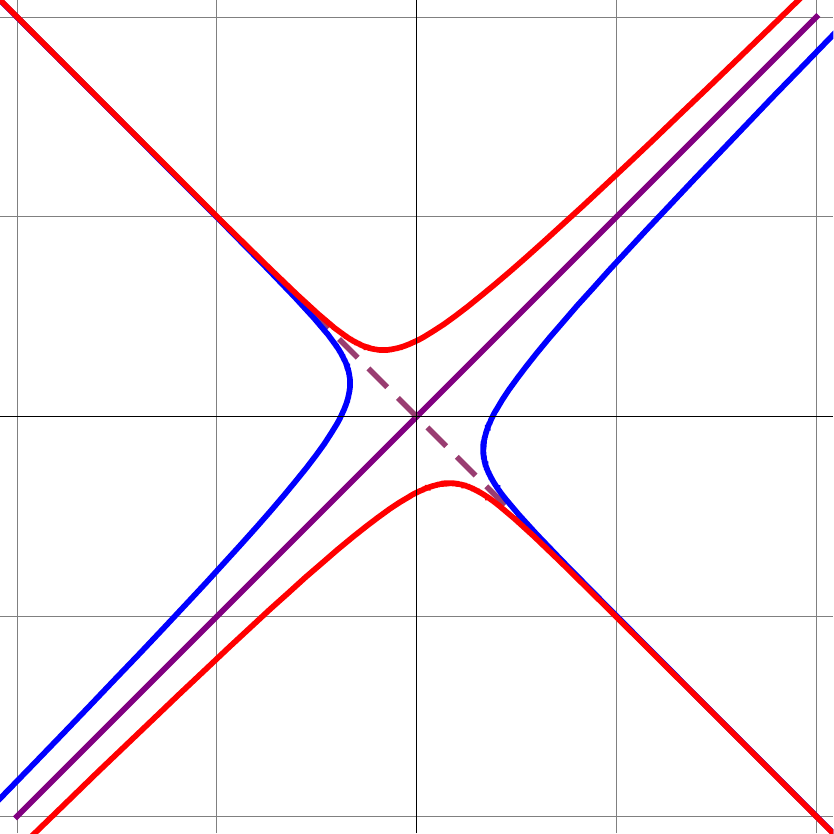}
\caption{Spacelike (blue) and timelike (red) elastic curve in $\lo^2$ with $\sigma = E= 0$.}
\label{Fig:elas0}
\end{center}
\end{figure}

In this case, $\mathcal K(v)=-\epsilon / v^2$ and we easily obtain $v(s)=1/s$, $u(s)=-\epsilon s^3/3$, $s\neq 0$. 
Using \eqref{cond_elastica}, its intrinsic equation is given by $\kappa (s) = 2/ s$, $s\neq 0$, that satisfies $4\dot \kappa^2=\kappa^4$. Hence, we get an elastic curve with $\sigma = E= 0$. See Figure \ref{Fig:elas0}.
\subsection{Case $c>0$}
\begin{figure}[h!]
\begin{center}
\includegraphics[height=4cm]{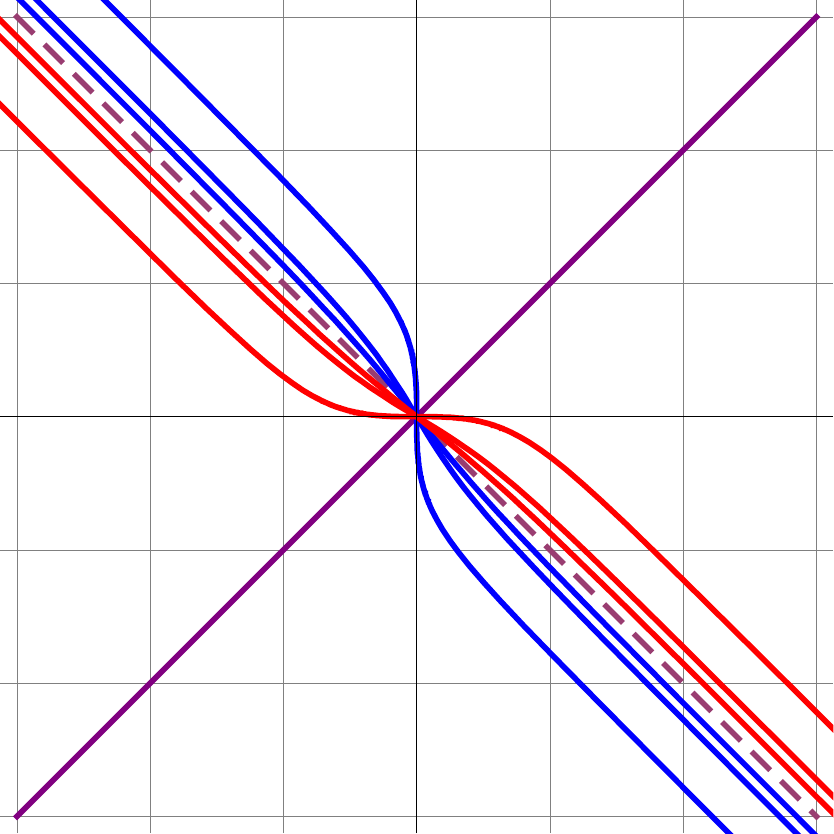}
\caption{Spacelike (blue) and timelike (red) elastic curves in $\lo^2$ with $\sigma = 4c>0$ and $E=4c^2$, $c=1,2,3$.}
\label{Fig:elasc-}
\end{center}
\end{figure}

From \eqref{s elas}, it is not difficult to get that $$v(s)=-\sqrt{c} \tan (\sqrt c \,s).$$ 
Using this in \eqref{u elas}, we deduce that $$u(s)=-\frac{\epsilon}{c} \left( \frac{s}{2}+\frac{\sin (2\sqrt c \,s)}{4 \sqrt{c}} \right).$$ 
From \eqref{cond_elastica} we obtain that the intrinsic equations of these curves are given by $\kappa (s)= -2\sqrt{c}\tan (\sqrt{c}s)$, $|s|< \pi/ 2 \sqrt{c}$.
It is a long exercise to check that it corresponds to an elastica under tension $\sigma = 4c>0$ and energy $E=4c^2$. See Figure \ref{Fig:elasc-}.
\subsection{Case $c<0$}
\begin{figure}[h!]
\begin{center}
\includegraphics[height=4cm]{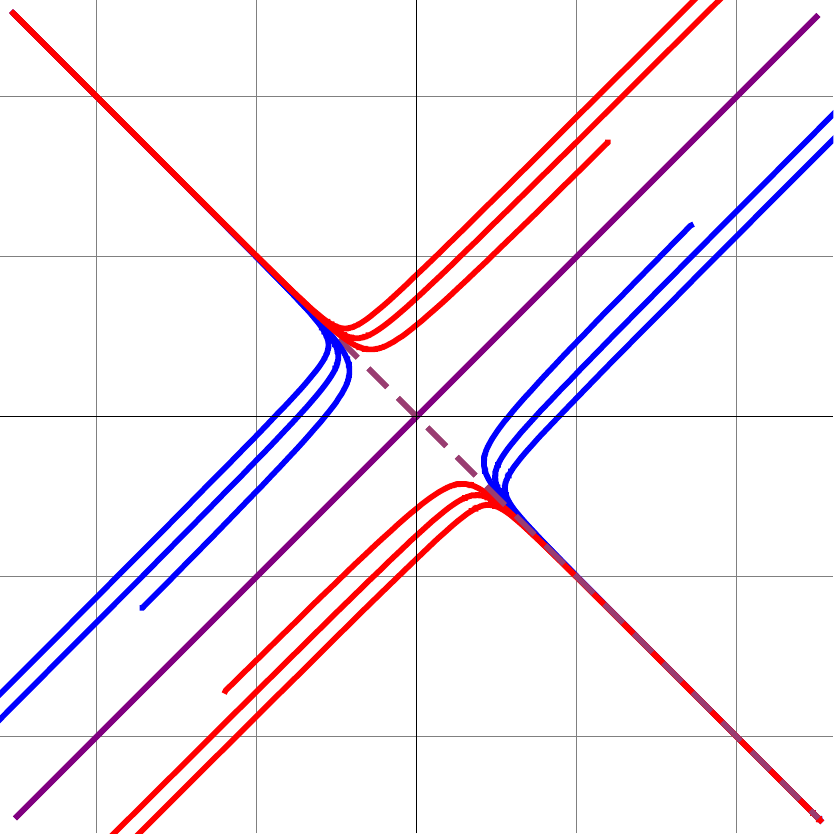}
\caption{Spacelike (blue) and timelike (red) elastic curves in $\lo^2$ with $\sigma = 4c<0$ and $E=4c^2$, $c=-1,-2,-3$.}
\label{Fig:elasc+}
\end{center}
\end{figure}

Using again \eqref{s elas}, we now obtain that $$v(s)=\sqrt{-c} \coth (\sqrt{-c} \,s).$$  Putting this in \eqref{u elas}, we deduce that $$u(s)=\frac{\epsilon}{c} \left(-\frac{s}{2}+\frac{\sinh (2\sqrt{-c} \,s)}{4 \sqrt{-c}} \right).$$
Taking into account \eqref{cond_elastica},  we have that the intrinsic equations of these curves are given by
$\kappa (s)= 2\sqrt{-c}\coth (\sqrt{-c}s)$, $s\neq 0$. It is again  a long straightforward computation to check that it corresponds to an elastica under tension $\sigma = 4c<0$ and energy $E=4c^2$. See Figure \ref{Fig:elasc+}.

In the last two cases, we remark that $E=\sigma ^2 /4$ and in this way we recover the special elastic curves described in Section 3 of \cite{CCIs18}.

\vspace{0.3cm}

\section{Curves in $\lo^2$ such that $\kappa (v)=a / v^2$, $a\neq 0$}

In this section we will study those spacelike and timelike curves in $\lo^2$ satisfying
\begin{equation}\label{cond_Ennepers}
\kappa(v)  = a / v^2,  a > 0.
\end{equation}
Given $\gamma =(u,v)$ satisfying \eqref{cond_Ennepers}, if we take $\hat \gamma = \frac{1}{a} (u,v)$ then, up to a dilation, we can only afford the condition
\begin{equation}\label{cond_Enn}
\kappa (v) = 1/ v^2,
\end{equation}
with $v\neq 0$.
Following Theorem~\ref{quadratures v}, we must consider the geometric linear momentum 
$$\mathcal K(v)=\frac{-\epsilon v}{c\, v-1}, \, c\in \R.$$

\subsection{Case $c=0$: $\mathcal K(v)=\epsilon v$}
\begin{figure}[h!]
\begin{center}
\includegraphics[height=4cm]{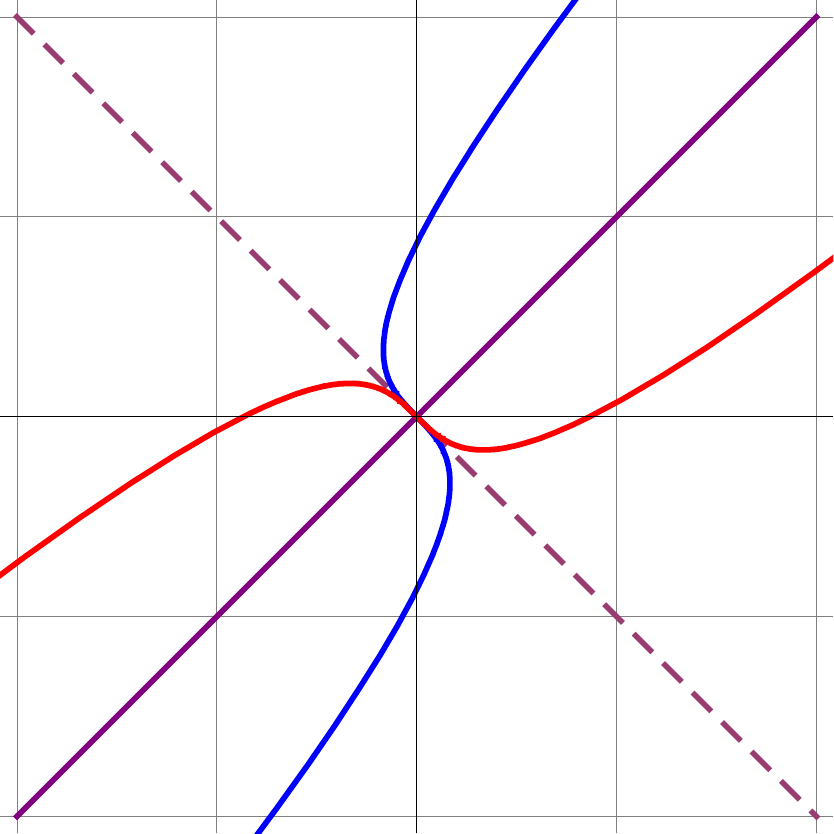}
\caption{Spacelike (blue) and timelike (red) curve in $\lo^2$ with $\mathcal K (v)=\epsilon v$, $\epsilon=\pm 1$.}
\label{Fig:Enn curves}
\end{center}
\end{figure}

We follow the steps described in Remark~\ref{c v} and we easily obtain that $v(s)=\sqrt{2s}$, $s>0$.
Recalling that $\mathcal K(v)=\epsilon v$, we get that $u(s)=2\epsilon\sqrt 2  s \sqrt s /3$, $s>0$.
We arrive at the graphs $u=\epsilon\, v^3 /3$, $v>0$, $\epsilon =\pm 1$. Using \eqref{cond_Enn}, their intrinsic equation is given by $\kappa (s)= \frac{1}{2s}, \, s>0$, for both of them (see Figure \ref{Fig:Enn curves}).

\begin{figure}[h!]
\begin{center}
\includegraphics[height=4cm]{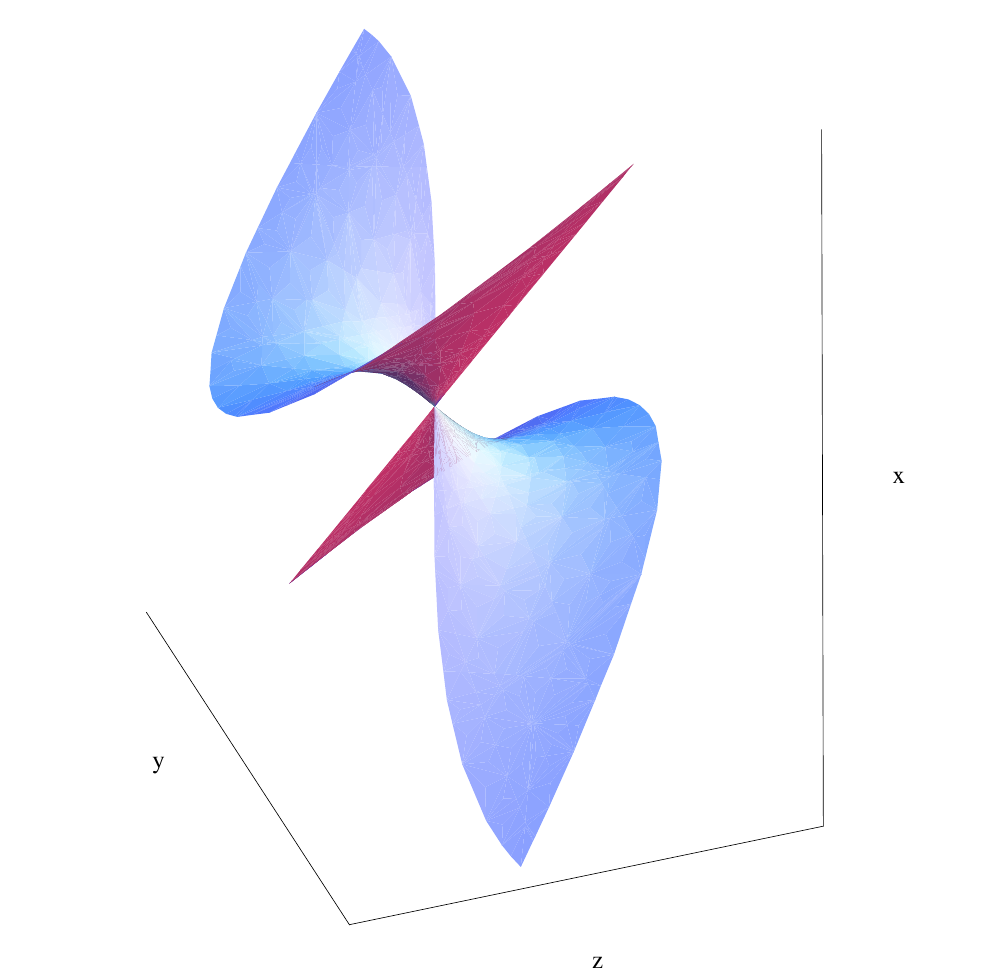}
\caption{Enneper's rotation surface of second kind in $\lo^3$.}
\label{Fig:Enn surface}
\end{center}
\end{figure}

On the other hand, Kobayashi introduced in \cite{Ko93}, by studying maximal rotation surfaces in $\lo^3:=(\R^3,-dx^2+dy^2+dz^2)$, the Enneper's surface of second kind. Concretely, 
Example 2.3 in \cite{Ko93} presents it (up to dilations) as a rotation surface with lightlike axis $(1,0,1)$ and generatrix curve $ x=\lambda(-t+t^3/3)$, $z=\lambda(t+t^3/3)$, $\lambda >0$, at the $xz$-plane (see Figure~\ref{Fig:Enn surface}).
Using \eqref{uv}, we notice that we have obtained exactly this curve (with $\lambda= 1/2$) in the spacelike case and 
we can conclude its following geometric characterization.
\begin{corollary}\label{cor:Ennep2}
The generatrix curve of the Enneper's surface of second kind, $u=v^3/3$, $v>0$, is the only spacelike curve (up to dilations and translations in the $u$-direction) in $\lo^2$ with geometric linear momentum $\mathcal K(v)=v$ (and curvature $\kappa (v)=1/v^2$). 
\end{corollary}
\subsection{Case $\mathcal K(v)=\frac{-\epsilon v}{c\, v-1}$, $c\neq 0$}
\begin{figure}[h!]
\begin{center}
\includegraphics[height=4cm]{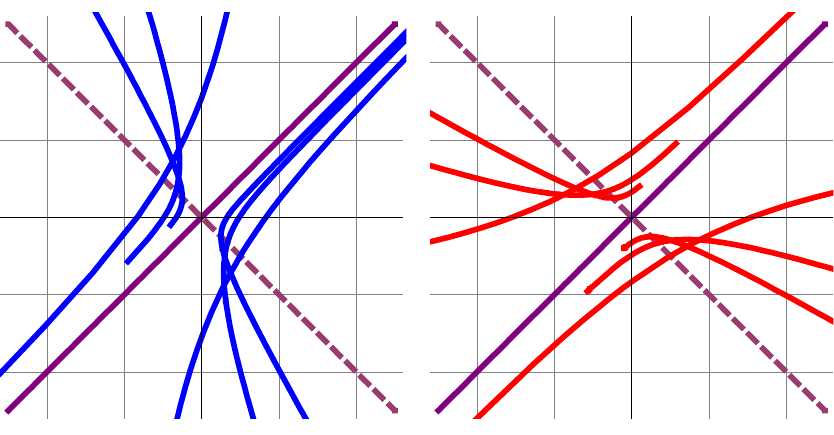}
\caption{Spacelike curves in $\lo^2$ with $\mathcal K (v)=-\frac{v}{c\, v-1}$ (left) and timelike curves in $\lo^2$ with $\mathcal K (v)=\frac{v}{c\, v-1}$ (right).}
\label{Fig:uplusminus}
\end{center}
\end{figure}

When $c\neq 0$, it is more difficult to get the arc parameter $s$ as a function of $v$; however, we can eliminate $ds$ using parts (ii) and (iii) in Remark \ref{c v}, obtaining after integration, the graphs
$$ u=u(v)=\frac{\epsilon}{c^3} \left(  c\, v -1 -\frac{1}{c\, v -1} +2  \log (c\, v -1)\right),$$
defined for $v>1/c$ if $c>0$ and for  $v<1/c$ if $c<0$ 
(see Figure \ref{Fig:uplusminus}).

\vspace{0.3cm}

\section{Curves in $\lo^2$ such that $\kappa (v)=a\, e^v $, $a\neq 0$}
In this section we will study those spacelike and timelike curves in $\lo^2$ satisfying
\begin{equation}\label{cond_exponen}
\kappa(v)  = a\,  e^v, a > 0. 
\end{equation}
Given $\gamma =(u,v)$ satisfying \eqref{cond_exponen}, if we take $\hat \gamma = (u,v+ \log a)$ then, up to a translation, we can only consider the condition
\begin{equation}\label{cond_exp}
\kappa (v) = e^v.
\end{equation}
Following Theorem~\ref{quadratures v}, we deal with the geometric linear momentum $$\mathcal K(v)=-\frac{\epsilon}{e^v+c}, \, c\in \R .$$

\subsection{Case $c=0$: $\mathcal K(v)\!=\!-\epsilon\, e^{-v}$.}
\begin{figure}[h!]
\begin{center}
\includegraphics[height=4cm]{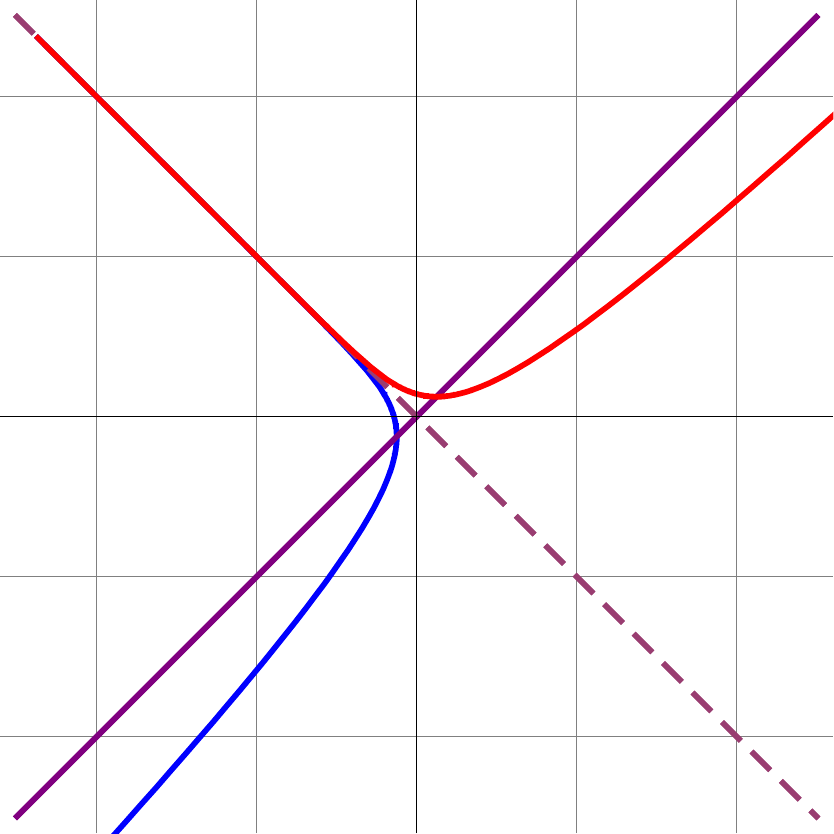}
\caption{Spacelike (blue) and timelike (red) curve in $\lo^2$ with $\mathcal K (v)=-\epsilon e^{-v}$, $\epsilon=\pm 1$.}
\label{Fig:Grims}
\end{center}
\end{figure}

Following the steps described in Remark~\ref{c v}, we easily obtain that $ v(s)=-\log s$, $s>0$.
And using that $\mathcal K(v)\!=\!-\epsilon \, e^{-v}$, we get that $u(s)=-\epsilon s^2/2$.
These curves correspond to the graphs $u=-\epsilon \, e^{-2v}/ 2$, $v\in \R$ (see Figure~\ref{Fig:Grims}). 
Using \eqref{cond_exp}, their intrinsic equations are given by $\kappa (s)= 1/s, \, s>0$.

It is straightforward to check that both curves satisfy the translating-type soliton equation $\kappa = g ((1,1),N)$.
Hence we have obtained in this section (see also Section 7.1 in \cite{CCI16}) certain Lorentzian versions of the grim-reaper curves of Euclidean plane. We will simply call them {\em Lorentzian grim-reapers}. As a summary, we conclude with the following geometric characterization of them.
\begin{corollary}\label{cor:grims}
The Lorentzian grim-reaper $u=- \, e^{-2v}/ 2$, $v\in \R$ (resp.\ $u= \, e^{-2v}/ 2$, $v\in \R$), is the only spacelike (resp.\ timelike)
curve ---up to translations in the $u$-direction--- in $\lo^2$ with geometric linear momentum $\mathcal K(v)=-  e^{-v}$
(resp.\  $\mathcal K(v)= e^{-v}$).
\end{corollary}

\subsection{Case $\mathcal K(v)\!=\!-\frac{\epsilon}{e^v+c}$, $c\neq 0$}
\begin{figure}[h!]
\begin{center}
\includegraphics[height=4cm]{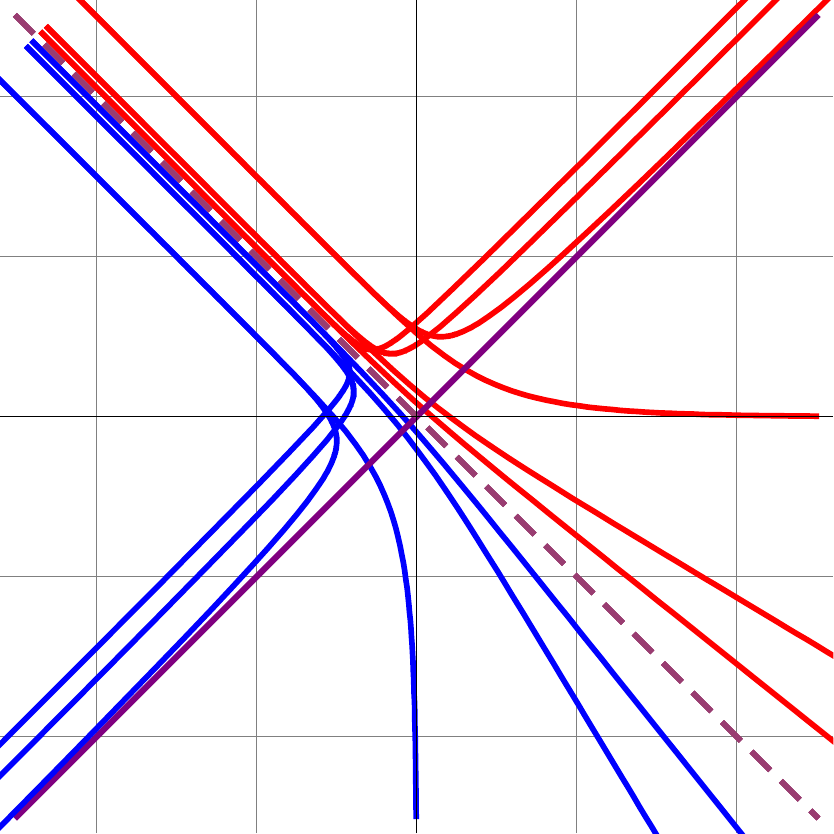}
\caption{Spacelike curves in $\lo^2$ (blue) and timelike curves (red) in $\lo^2$ with $\mathcal K (v)=-\frac{\epsilon}{e^v+c}$, $c\neq 0$, $\epsilon=\pm 1$.}
\label{Fig:exp v}
\end{center}
\end{figure}

When $c\neq 0$, a straightforward computation, after solving the corresponding integrations, leads to the following curves:
$$ v(s)=\log \frac{c}{e^{cs}-1}, \ u(s)=-\frac{\epsilon}{c}\left(  s+\frac{1}{c\, e^{cs}} \right), \, s>0. $$
Using \eqref{cond_exp}, we deduce that their intrinsic equations are given by 
$  \kappa (s)= \frac{c}{e^{cs}-1}$, $s>0$, $c\neq 0$  (see Figure~\ref{Fig:exp v}).

\end{document}